%
%
%
%
%
%
\magnification=\magstep1
\input amstex
\documentstyle{amsppt}
\voffset-3pc

\define\C{{\Bbb C}}

\define\dee{\partial}
\redefine\O{\Omega}
\redefine\phi{\varphi}
\define\Obar{\overline{\Omega}}
\define\Ot{\widetilde\Omega}
\define\Oh{\widehat\Omega}

\NoRunningHeads
\topmatter
\title
Complexity in complex analysis
\endtitle
\author Steven R. Bell${}^*$ \endauthor
\thanks ${}^*$Research supported by NSF grant DMS-9623098 \endthanks
\address
Mathematics Department, Purdue University, West Lafayette, IN  47907 USA
\endaddress
\email bell\@math.purdue.edu \endemail
\abstract
We show that the classical kernel and domain functions associated to an
$n$-connected domain in the plane are all given by rational combinations of
three or {\it fewer\/} holomorphic functions of {\it one\/} complex
variable.  We characterize those domains for which the classical
functions are given by rational combinations of only two or fewer functions
of one complex variable.  Such domains turn out to have the property that
their classical domain functions all extend to be meromorphic functions on a
compact Riemann surface, and this condition will be shown to be
equivalent to the condition that an Ahlfors map and its derivative are
algebraically dependent.  We also show how many of these results can
be generalized to finite Riemann surfaces.
\endabstract
\endtopmatter
\document

\hyphenation{bi-hol-o-mor-phic}
\hyphenation{hol-o-mor-phic}

\subhead 1. Introduction \endsubhead
On a simply connected domain $\O\ne\C$ in the plane,
the classical Bergman kernel $K(z,w)$ associated to $\O$ is given by
$$K(z,w)=
\frac{f_a'(z)\overline{f_a'(w)}}
{\pi (1-f_a(z)\overline{f_a(w)})^2},$$
where $f_a(z)$ is the Riemann mapping function mapping $\O$ one-to-one onto
the unit disc $D_1(0)$ with $f_a(a)=0$ and $f_a'(a)>0$.  Thus, the Bergman
kernel is a rational combination of just {\it two\/} holomorphic functions
of {\it one\/} complex variable.  I have recently proved in \cite{9} that
the Bergman kernel and many other objects of potential theory associated
to a finitely multiply connected domain are rational combinations
of only {\it three\/} holomorphic function of {\it one\/} complex variable,
namely two Ahlfors maps plus the derivative of one Ahlfors map.  In
\cite{8}, I proved that the three functions of $z$ given by $S(z,A_j)$,
where $S(z,w)$ is the Szeg\H o kernel and $A_j$, $j=1,2,3$, are three
fixed points in the domain, generate the Szeg\H o kernel, the Bergman
kernel, and many other objects of potential theory.  In
this paper, I shall unify these results and I shall show that the three
functions of one variable that generate these classical functions can be
taken from a rather long list of functions.  I also showed in \cite{10}
that there exist certain multiply connected domains in the plane that
are particularly simple in the sense that their Bergman and Szeg\H o
kernels are generated by only {\it two\/} functions
of one complex variable.  In this paper, I shall present all these results
in a single improved framework and I will thereby be able to characterize
those domains whose kernel functions are particularly simple.  These domains
are characterized by a condition on the Ahlfors map which turns
out to be equivalent to the condition that the Bergman kernel and other
objects of potential theory associated to the domain extend to a compact
Riemann surface as single valued meromorphic functions.  At the heart
of these results is a relationship between the Szeg\H o kernel of a
domain and proper holomorphic maps of the domain to the unit disc.
Because we can prove similar relationships on any finite Riemann surface,
we will see that many of these results can be generalized to this setting.

\subhead 2. Statement of main results \endsubhead
Because the results of this paper are easiest to state and prove
for domains in the plane and because all of the main technical
advances can be understood in this setting, we shall devote the
first part of the paper to planar domains.  In the last section of the
paper, we sketch the rather routine argument to generalize the
results to the case of finite Riemann surfaces.

We shall be able to prove our main results for finitely connected
domains in the plane such that no boundary component is a point.  For
the moment, however, assume that $\O$ is a bounded $n$-connected domain in
the plane bounded by $n$ non-intersecting $C^\infty$ smooth simple
closed curves, $\gamma_j$, $j=1,...n$.  For a point $a$ in
$\O$, let $f_a$ denote the Ahlfors map associated to the pair ($\O,a)$.
This map is an $n$-to-one (counting multiplicities) proper holomorphic
mapping of $\O$ onto the unit disc.  Furthermore, $f_a(a)=0$, and $f_a$
is the unique holomorphic function mapping $\O$ into the unit disc
maximizing the quantity $|f_a'(a)|$ with $f_a'(a)>0$.  It is well
known that $f_a$ extends $C^\infty$ smoothly up to the boundary $b\O$
of $\O$ and that $|f_a|=1$ on $b\O$.

Let $K(z,w)$ denote the Bergman kernel associated to $\O$ and
let $\omega_j$ denote the harmonic measure function which is
harmonic on $\O$ with boundary values of one on the boundary
component $\gamma_j$ and zero on the other boundary components.
Let $F_j'(z)$ denote the holomorphic function given by
$(1/2)(\dee/\dee z)\omega_j(z)$.  (The prime is traditional;
$F_j'$ is not the derivative of a holomorphic function on $\O$.)
Let $G(z,a)$ denote the classical Green's function associated to
$\O$.

The Szeg\H o kernel associated to $\O$ is the kernel for the
orthogonal projection of $L^2(b\O)$ (with respect to arc length
measure) onto the subspace consisting of $L^2$ boundary values
of holomorphic functions on $\O$, i.e., the Hardy space $H^2(b\O)$.

When $\O$ is merely an $n$-connected domain in the plane such that
no boundary component is a point, we define the Szeg\H o kernel associated
to $\O$ as follows.  There exists a biholomorphic mapping $\Phi$ mapping
$\O$ one-to-one onto a bounded domain $\O^a$ in the plane with
real analytic boundary.  The standard construction yields a domain
$\O^a$ that is a bounded $n$-connected domain with $C^\infty$ smooth
real analytic boundary whose boundary consists of $n$ non-intersecting
simple closed real analytic curves.  The function $\Phi'$ has a
single valued holomorphic square root on $\O$ (see \cite{4, page 43}).
Let superscript $a$'s indicate that a kernel function is associated
to $\O^a$.  Kernels without superscripts are associated to $\O$.
The Szeg\H o kernel associated to $\O$ is {\it defined\/} via
$$S(z,w)=\sqrt{\Phi'(z)}\
S^a(\Phi(z),\Phi(w))\overline{\sqrt{\Phi'(w)}}.$$
We shall also define the Garabedian kernel associated to $\O$ via the
natural transformation formula,
$$L(z,w)=\sqrt{\Phi'(z)}\ L^a(\Phi(z),\Phi(w))\sqrt{\Phi'(w)}.$$
(See \cite{4, page~24} or \S3 of this paper for the definition of
the Garabedian kernel in a smoothly bounded domain.)  Various other
transformation formulas hold for the objects of potential theory
mentioned above.  For example, the Bergman kernels transform via
$$K(z,w)=\Phi'(z)K^a(\Phi(z),\Phi(w))\overline{\Phi'(w)},$$
and the Green's functions satisfy
$$G(z,w)=G^a(\Phi(z),\Phi(w)),$$
and the functions associated to harmonic measure satisfy
$$\omega_j(z)=\omega_j^a(\Phi(z))\quad\text{ and }\quad
F_j'(z)=\Phi'(z) {F_j^a}'(\Phi(z))$$
(provided, of course, that we stipulate that the boundary components
have been numbered so that $\Phi$ maps the $j$-th boundary component of
$\O$ to the $j$-th boundary component of $\O^a$).
The Ahlfors map $f_b$ associated to a point $b\in\O$ is the holomorphic
function mapping $\O$ into $D_1(0)$ with $|f_b'(b)|$ maximal and
$f_b'(b)>0$.  It is easy to see that the Ahlfors map satisfies
$$f_b(z)=\lambda f_{\Phi(b)}^a(\Phi(z))$$
for some unimodular constant $\lambda$ and it follows that
$f_b(z)$ is a proper holomorphic mapping of $\O$ onto $D_1(0)$.
The double of $\O$ may also be defined in a standard way by using
$\Phi$ and the double of $\O^a$.

The work in this paper is motivated by the following result from \cite{9}.

\proclaim{Theorem 2.1}
Suppose $\O$ is an $n$-connected domain in the plane such that no boundary
component is a point.  There exist points $a$ and $b$ in $\O$ and complex
rational functions $R$ and $Q$ of four complex variables such that the
Bergman kernel can be expressed in terms of the two Ahlfors maps $f_a$
and $f_b$ via
$$K(z,w)=f_a'(z)\overline{f_a'(w)}
R(f_a(z),f_b(z),\overline{f_a(w)},\overline{f_b(w)})$$
and the Szeg\H o kernel can be expressed via
$$S(z,w)^2=f_a'(z)\overline{f_a'(w)}
Q(f_a(z),f_b(z),\overline{f_a(w)},\overline{f_b(w)}).$$
Alternatively, the Szeg\H o kernel can be expressed via
$$S(z,w)=S(z,a)S(a,w)
Q_2(f_a(z),f_b(z),\overline{f_a(w)},\overline{f_b(w)})$$
where $Q_2$ is rational.  Furthermore, the functions $F_j'$ are given by
$$F_j'(z)=f_a'(z)P(f_a(z),f_b(z))$$
where $P$ is a rational function of two complex variables.  Also,
every proper holomorphic mapping of $\O$ onto the unit disc is a
rational combination of $f_a$ and $f_b$.
\endproclaim

There are similar formulas for the complex derivative of the
Green's function $(\dee/\dee z)G(w,z)$ and for the Poisson kernel
given in \cite{9}.

These results show that the analytic objects of potential theory
and conformal mapping in a multiply connected domain are all
given as rational combinations of only {\it three\/} holomorphic functions
of one complex variable.  In \cite{10}, I showed that, when the
kernel functions are algebraic, these analytic objects are given
by rational combinations of only {\it two\/} holomorphic functions of
one complex variable.  In the present paper, we characterize this
property in the following theorem.  We shall say that a kernel
function such as the Bergman kernel $K(z,w)$ is given as a
rational combination of two functions of one complex variable
on $\O$ if there exist holomorphic functions $G_1$ and $G_2$ on
$\O$ such that $K(z,w)$ can be written as a rational combination
of $G_1(z)$, $G_2(z)$, $\overline{G_1(w)}$, and $\overline{G_2(w)}$.

\proclaim{Theorem 2.2}
Suppose $\O$ is an $n$-connected domain in the plane with $n>1$ such
that no boundary component is a point.  The following conditions are
equivalent.
\roster
\item The Bergman kernel $K(z,w)$ is given as a rational combination
of only two functions of one complex variable on $\O$.
\item The Szeg\H o kernel $S(z,w)$ is given as a rational combination
of only two functions of one complex variable on $\O$.
\item The Bergman kernel, Szeg\H o kernel, the classical functions
$F_j'(z)$ associated to $\O$, and all proper holomorphic mappings of
$\O$ onto the unit disc are all given as rational combinations of the
same two functions of one complex variable on $\O$.
\item The domain $\O$ can be realized as a subdomain of a compact
Riemann surface ${\Cal R}$ such that $S(z,w)$ and $K(z,w)$ extend
to ${\Cal R}\times{\Cal R}$ as single valued
meromorphic functions, and as such, can be expressed as rational
combinations of any two functions on $\O$ which extend to $\Cal R$
to form a primitive pair for $\Cal R$.  Also, every proper holomorphic
mapping of $\O$ onto the unit disc as well as the functions $F_k'(z)$,
$k=1,\dots,n-1$, extend to be single valued meromorphic functions on
${\Cal R}$.  Furthermore, the complement of $\O$ in ${\Cal R}$ is
connected.
\item There exists a single proper holomorphic mapping of $\O$ onto
the unit disc which satisfies an identity of the form
$$P(f'(z), f(z))=0$$
on $\O$ for some complex polynomial $P$, i.e., $f$ and $f'$ are
algebraically dependent.
\item For every proper holomorphic mapping $f$ of $\O$ onto
the unit disc, $f$ and $f'$ are algebraically dependent.
\endroster
\endproclaim

A novel feature of the proof of Theorem~2.2 is the use of the
Zariski-Castelnuovo Theorem from Algebra.

We remark that it is always possible to choose the two functions
that form the primitive pair in condition (4) to be holomorphic
on $\O$ (see Farkas and Kra \cite{15}).

It shall fall out from the proof of Theorem~2.2 that if the Bergman
or Szeg\H o kernel associated to a finitely connected domain with no
pointlike boundary components is generated by a {\it single\/} function
of one variable, then the domain must be simply connected.  In fact,
if the domain under study is simply connected, the proof of Theorem~2.2
becomes considerably easier and a similar result can be established
without using L\"uroth's Theorem in place of the Zariski-Castelnuovo Theorem.

\proclaim{Theorem 2.2a}
Suppose $\O$ is a simply connected domain in the plane not equal to
$\C$.  The following conditions are equivalent.
\roster
\item The Bergman or Szeg\H o kernels are given as rational combinations
of only one function of one complex variable on $\O$.
\item The domain $\O$ can be realized as a subdomain of a compact
Riemann surface ${\Cal R}$ such that $S(z,w)$ and $K(z,w)$ extend
to ${\Cal R}\times{\Cal R}$ as single valued
meromorphic functions.
\item A Riemann mapping of $\O$ onto
the unit disc and its derivative are algebraically dependent.
\endroster
\endproclaim

Many frequently encountered domains satisfy the polynomial condition
(5) in Theorem~2.2 for some proper map.  For example, this condition
holds for generalized quadrature domains of the type studied by Aharonov
and Shapiro (see \cite{1, page~64} and \cite{18}) and it holds for
domains that have algebraic kernel functions (see \cite{{7,10}}).
This last class of domains is the same as the set of domains given by
a connected component $\O$ of a set of the form $\{z:|f(z)|<1\}$ where
$f$ is an algebraic function without singularities in $\O$, i.e., $\O$
is a domain such there exists an algebraic proper holomorphic mapping
from it onto the unit disc (see \cite{7}).

Theorem~2.1 is actually a special case of a more general theorem
that we prove here.  Let $\Oh$ denote the double of $\O$ and let $R(z)$ denote
the antiholomorphic involution on $\Oh$ which fixes the boundary of
$\O$.  Let $\Ot=R(\O)$ denote the reflection of $\O$ in $\Oh$ across
the boundary.  The double of $\O$ is a compact Riemann surface and
hence the field of meromorphic functions on $\Oh$ is generated by a
pair of meromorphic functions $(G_1,G_2)$ on $\Oh$ known as a primitive
pair (see Farkas and Kra \cite{15}).  The construction of primitive
pairs given in \cite{15} shows that we may assume that $G_1$ and $G_2$
are holomorphic on $\O\subset\Oh$.

We now define a class~$\Cal A$ of meromorphic functions on $\O$.  Recall
that $G(z,a)$ denotes the classical Green's function associated to $\O$.

The class~$\Cal A$ consists of

\roster
\item the functions $F_j'(z)$, $j=1,\dots,n$,
\item functions of $z$ of the form $\frac{\dee}{\dee z}G(z,a)$ for fixed
$a$ in $\O$,
\item functions of $z$ of the form $D_a\frac{\dee}{\dee z}G(z,a)$ where
$D_a$ denotes a differential operator of the form $\frac{\dee^n}
{\dee a^n}$ or $\frac{\dee^n}{\dee\bar a^n}$, and $a$ is a fixed point
in $\O$,
\item functions of $z$ of the form $S(z,a_1)S(z,a_2)$ where $a_1$ and
$a_2$ are fixed points in $\O$,
\item and linear combinations of functions above.
\endroster

If $\O$ has $C^\infty$ smooth boundary, we allow the points $a$ in
(2) and (3) of the definition of the class $\Cal A$ to be in the larger
set $\Obar$.

Theorem~2.1 will be generalized as follows.

\proclaim{Theorem 2.3}
Suppose that $\O$ is a finitely connected domain in the plane such
that no boundary component is a point.  Let $G_1$ and $G_2$ denote any
two meromorphic functions on $\O$ that extend to the double of $\O$ to
form a primitive pair, and let $A(z)$ denote any function from the class~
$\Cal A$ other than the zero function.  The Bergman kernel associated
to $\O$ can be expressed as
$$K(z,w)=A(z)\overline{A(w)}
R_1(G_1(z),G_2(z),\overline{G_1(w)},\overline{G_2(w)})$$
where $R_1$ is a complex rational function of four complex variables.
Similarly, the Szeg\H o kernel can be expressed as
$$S(z,w)^2= A(z)\overline{A(w)}
R_2(G_1(z),G_2(z),\overline{G_1(w)},\overline{G_2(w)})$$
where $R_2$ is rational, and
the functions $F_j'$ can be expressed
$$F_j'(z)=A(z)R_3(G_1(z),G_2(z))$$
where $R_3$ is rational.  Furthermore, every proper holomorphic
mapping of $\O$ onto the unit disc is a rational combination of
$G_1$ and $G_2$.
\endproclaim

The proof of Theorem~2.3 hinges on the fact proved in \S6 that, when
$\O$ has real analytic boundary, functions in the class~$\Cal A$ can
be seen to be equal to meromorphic functions $H$ on $\Obar$ that satisfy
an identity of the form
$$H(z)T(z)=\overline{J(z)T(z)}$$
for $z\in b\O$ where $J(z)$ is another meromorphic function on $\Obar$
and $T(z)$ represents the complex unit tangent vector at $z$ pointing
in the direction of the standard orientation of $b\O$.  We prove in
\S6 that the class~$\Cal A$ is the largest set of functions with this
property.

It is remarkable how many different common functions fall into the
class~$\Cal A$.  For example, the Bergman kernel is in $\Cal A$ because
it is related to the classical Green's function via
(\cite{14, page~62}, see also \cite{4, page~131})
$$K(z,w)=-\frac{2}{\pi}\frac{\dee^2 G(z,w)}{\dee z\dee\bar w}.$$
Another kernel function on $\O\times\O$ that we shall need is given by
$$\Lambda(z,w)=-\frac{2}{\pi}\frac{\dee^2 G(z,w)}{\dee z\dee w}.$$
(In the literature, this function is sometimes written as $L(z,w)$ with
anywhere between zero and three tildes and/or hats over the top.
We have chosen the symbol $\Lambda$ here to avoid
confusion with our notation for the Garabedian kernel above.)
Note that $\Lambda(z,w)$ is in $\Cal A$ as a function of $z$ for each
fixed $w$ in $\O$.  Furthermore, given a proper holomorphic map $f$
from $\O$ onto the unit disc, the quotient $f'/f$ is in $\Cal A$
by virtue of the fact that $\ln|f|$ is a linear combination of
Green's functions.  Since $\bar f=1/f$ on $b\O$, we shall be able to
see that $f'(z)T(z)=-\overline{T(z)f'(z)/f(z)^2}$ for $z\in b\O$ on
a smooth domain, and we will be able to deduce that $f'$ itself is in
$\Cal A$ from Lemma~6.3 below.

It is also remarkable the variety of functions that can appear
as members of a primitive pair.  The argument in \cite{9, page~332}
reveals that for any proper holomorphic map $f_1$ from $\O$ onto the unit
disc, there exists a second proper map $f_2:\O\to D_1(0)$ such
that $(f_1,f_2)$ forms a primitive pair for the double of $\O$.  (In
fact, the second map can be taken to be an Ahlfors map.)  When $f_1'/f_1$
is taken to be equal to $A(z)$ in Theorem~2.3 and $G_1$ and $G_2$
are taken to be this primitive pair, the Bergman kernel $K(z,w)$ is
expressed as $f_1'(z)\overline{f_1'(w)}$ times a rational function of
$f_1(z)$, $f_2(z)$, $\overline{f_1(w)}$ and $\overline{f_2(w)}$.  This
result makes it very easy to prove
the result given in \cite{7} that the Bergman kernel associated to a
finitely connected domain is an algebraic function if and only if there
exists a proper holomorphic map from the domain onto the unit disc
which is algebraic.  Indeed, if $f_1$ is algebraic, then $f_1'$ is
algebraic, and since $f_1$ and $f_2$ extend to the double to form
a primitive pair, they must be algebraically dependent (see Farkas
and Kra \cite{15, page~248}).  Hence, $f_2$ is an algebraic function
of $f_1$ and this shows that $K(z,w)$ is algebraic.  The reverse
implication is even easier (see \cite{{7,11,12}}).  Similarly, the
statement about the Szeg\H o kernel in Theorem~2.3 yields an easy proof of
the result given in \cite{7} that the Szeg\H o kernel associated to a
finitely connected domain is an algebraic function if and only if there
exists a proper holomorphic map from the domain onto the unit disc
which is algebraic.

It is shown in \cite{8}
that, for almost any three points $a_1$, $a_2$, and $b$ in $\O$,
the two quotients $K(z,a_1)/K(z,b)$ and $K(z,a_2)/K(z,b)$ extend
to the double of $\O$ and form a primitive pair.  When these two functions
are used as $G_1$ and $G_2$ and when $K(z,b)$ is used as the function
$A(z)$ from the class~$\Cal A$, we find that $K(z,w)$ is a rational combination
of $K(z,a_1)$, $K(z,a_2)$, and $K(z,b)$ and conjugates of $K(w,a_1)$,
$K(w,a_2)$, and $K(w,b)$ (as was also shown in \cite{8}).  It now follows
from Theorem~2.3 that all the functions in the class~$\Cal A$ are rational
combinations of $K(z,a_1)$, $K(z,a_2)$ and $K(z,b)$.  This reinforces the
idea that the Bergman kernel contains {\it everything\/} there is
to know about a domain.

Let $G_z(z,w)$ denote the derivative $(\dee/\dee z)G(z,w)$ of the Green's
function.  Similar reasoning to that given in \cite{8} shows
that, for almost any three points $a_1$, $a_2$, and $b$ in $\O$, the two
quotients $G_z(z,a_1)/G_z(z,b)$ and $G_z(z,a_2)/G_z(z,b)$ form a primitive
pair for the double of $\O$.  Hence, Theorem~2.3 yields that the
Bergman kernel is a rational combination of the three functions of
one variable given by $G_z(z,a_1)$, $G_z(z,a_2)$, and $G_z(z,b)$.
A similar result has been proved for the Szeg\H o kernel in \cite{8}.
Let us summarize these results in the statement of the following theorem.

\proclaim{Theorem 2.4}
Suppose that $\O$ is a finitely connected domain
in the plane such that no boundary component is a point.
The Bergman kernel $K(z,w)$ associated to $\O$ can be
expressed as a rational combination of the three functions
\roster
\item
$K(z,A_1)$, $K(z,A_2)$, $K(z,A_3)$, or
\item
$S(z,A_1)$, $S(z,A_2)$, $S(z,A_3)$, or
\item
$G_z(z,A_1)$, $G_z(z,A_2)$, $G_z(z,A_3)$,
\endroster
where $A_1$, $A_2$, and $A_3$ are three fixed points in $\O$.
\endproclaim

Recall that $L(z,a)$ denotes the Garabedian kernel associated to $\O$.
We now define another class of meromorphic functions on $\O$ relevant
to the Szeg\H o kernel.  The class~$\Cal B$ consists of

\roster
\item functions of $z$ of the form $S(z,a)$ or $L(z,a)$ for fixed points
$a$ in $\O$,
\item functions of $z$ of the form $\frac{\dee^m}{\dee\bar a^m }S(z,a)$
or $\frac{\dee^m}{\dee a^m }L(z,a)$ for fixed points $a$ in $\O$,
\item and linear combinations of functions above.
\endroster

If $\O$ has $C^\infty$ smooth boundary, we allow the points $a$ in
(1) and (2) in the definition of the class $\Cal B$ to be in the larger
set $\Obar$.

The Szeg\H o kernel can be expressed in a manner similar to the Bergman
kernel as follows.

\proclaim{Theorem 2.5}
Suppose that $\O$ is a finitely connected domain in the plane such
that no boundary component is a point.  Let $G_1$ and $G_2$ denote any
two meromorphic functions on $\O$ that extend to the double of $\O$ to
form a primitive pair, and let $B(z)$ denote any function from the
class~$\Cal B$ other than the zero function.  The Szeg\H o kernel
associated to $\O$ can be expressed as
$$S(z,w)=B(z)\overline{B(w)}
R(G_1(z),G_2(z),\overline{G_1(w)},\overline{G_2(w)})$$
where $R$ is a complex rational function of four complex variables.
\endproclaim

The proof of Theorem~2.5 hinges on the fact proved in \S6 that, when
$\O$ has real analytic boundary, functions in the class~$\Cal B$ can
be seen to be equal to meromorphic functions $G$ on $\Obar$ that satisfy
an identity of the form
$$G(z)=\overline{H(z)T(z)}$$
for $z\in b\O$ where $H(z)$
is another meromorphic function on $\Obar$ and, as before, $T(z)$
represents the complex unit tangent vector at $z$ pointing
in the direction of the standard orientation of $b\O$.  We prove in
\S6 that the class~$\Cal B$ is the largest set of functions with this
property.

It is shown in \cite{8} that, for almost any three
points $a_1$, $a_2$, and $b$ in $\O$, the two quotients $S(z,a_1)/S(z,b)$
and $S(z,a_2)/S(z,b)$ form a primitive pair for the double of $\O$.  When
these two functions are used as $G_1$ and $G_2$ and when $S(z,b)$ is used
as $B(z)$ in Theorem~2.5, we find that $S(z,w)$ is a rational combination
of $S(z,a_1)$, $S(z,a_2)$ and $S(z,b)$, and the conjugates of $S(w,a_1)$,
$S(w,a_2)$ and $S(w,b)$, as stated in Theorem~2.4.  It now also follows
from Theorem~2.5 that all the
functions in the class~$\Cal B$ are rational combinations of $S(z,a_1)$,
$S(z,a_2)$ and $S(z,b)$.  It can be shown that the class $\Cal A$ is
equal to the complex linear span of the set of products of two functions
in the class~$\Cal B$.   Hence, it is reasonable to say that
the Szeg\H o kernel contains {\it even more\/} information about
a domain than the Bergman kernel does.

It is reasonable to wonder why on earth one might want to take complicated
linear combinations for $A(z)$ or $B(z)$ in Theorems~2.3 and 2.5 instead of a
single simple function from the classes.  However, if $\O$ is a
quadrature domain in the sense that
$$\iint_\O h(z)\ dA = \sum_{j=1}^N h(w_j)$$
for finitely many fixed points $w_j$ in $\O$ and all holomorphic functions
$h$ in the Bergman space, then
$$\sum_{j=1}^N K(z,w_j)\equiv1$$
and we may take $A(z)\equiv1$ in Theorem~2.3.  Hence, if $\O$ is a
quadrature domain in this sense, then the Bergman kernel is a rational
combination of any two functions that extend to the double of $\O$ to
form a primitive pair.  Similarly, if $\O$ is a
quadrature domain in the sense that
$$\int_{b\O} h(z)\ ds = \sum_{j=1}^N h(w_j)$$
for finitely many fixed points $w_j$ in $\O$ and all holomorphic functions
$h$ in the Hardy space, then
$$\sum_{j=1}^N S(z,w_j)\equiv1$$
and we may take $B(z)\equiv1$ in Theorem~2.5.  Hence, if $\O$ is a
quadrature domain in this sense, then the Szeg\H o {\it and\/} Bergman kernels
are rational combinations of any two functions that extend to the double
of $\O$ to form a primitive pair.

Most of the results of this paper depend on a special formula for the Szeg\H o
kernel.  It is a standard construction to produce the Szeg\H o projection
and kernel with respect to a weight function on the boundary of
$\O$ when the boundary of $\O$ is sufficiently smooth.  Suppose that
$\O$ is a domain in the plane bounded by finitely many $C^\infty$
smooth curves.  (It will be clear from the proofs given later that
this smoothness condition could be greatly relaxed, but we will be
able to reap enough consequences in the smooth case that we will
not bother trying to generalize the result here.)  Given a positive
$C^\infty$ weight function $\phi$ on $b\O$, let $S(z,w)$ denote the Szeg\H o
kernel function defined on $\O\times b\O$ which reproduces holomorphic
functions on $\O$ with respect to the weight $\phi$ in the sense that
$$h(z)=\int_{w\in b\O}S(z,w)h(w)\,\phi(w)\,ds$$
for points $z$ in $\O$ and holomorphic functions $h$ on $\O$ that are
in the Hardy space associated to $\O$.  We shall prove that the weighted
Szeg\H o kernel $S(z,w)$ is a rational combination of finitely many
holomorphic functions of one complex variable on $\O$ in the following
precise sense.

Let $S_0(z,w)$ denote $S(z,w)$, and let
$S_{\bar n}(z,w)$ denote $(\dee^n/\dee\bar w^n)S(z,w)$.

\proclaim{Theorem 2.6}
Suppose that $\O$ is a domain in the plane bounded by finitely many
$C^\infty$ smooth curves.  Suppose that $f:\O\to D_1(0)$ is a
proper holomorphic mapping that has zeroes at $a_1,\dots,a_N$ with
multiplicities $M(1),\dots,M(N)$, respectively.  Given a positive
$C^\infty$ weight function $\phi$ on $b\O$, the weighted Szeg\H o kernel
$S(z,w)$ with respect to $\phi$ satisfies
$$S(z,w)=\frac{1}{1-f(z)\overline{f(w)}}\left(
\sum_{i,j=1}^{N}
\sum_{n=0}^{M(i)}
\sum_{m=0}^{M(j)}
c_{ijnm}S_{\bar n}(z,a_i)\,\overline{S_{\bar
m}(w,a_j)}\right)$$
for some constants $c_{ijnm}$.
\endproclaim

The constants $c_{ijnm}$ will be described more fully in the
proof of this theorem given in \S8.

When the weight $\phi$ is taken to be the Poisson kernel associated to
a point $a$ in $\O$, we will be able to say more about the associated
Szeg\H o kernel and how it is related to the objects of potential
theory, see \S10.  This viewpoint will also allow us to find
interesting generalizations to the case where $\O$ is a finite Riemann
surface, see \S11.

\subhead 3. Preliminaries \endsubhead
For the moment, suppose that $\O$ is a bounded $n$-connected domain in the
plane with $C^\infty$ smooth boundary, i.e., a domain whose boundary
$b\O$ is given by finitely many non-intersecting $C^\infty$ simple closed
curves.  Let $\gamma_j$, $j=1,\dots,n$, denote the $n$ non-intersecting
$C^\infty$ simple closed curves which define the boundary of $\O$, and
suppose that $\gamma_j$ is parameterized in the standard sense by $z_j(t)$.
Let $T(z)$ be the $C^\infty$ function
defined on $b\O$ such that $T(z)$ is the complex number representing
the unit tangent vector at $z\in b\O$ pointing in the direction of
the standard orientation.  This complex unit tangent vector function
is characterized by the equation $T(z_j(t))=z_j'(t)/|z_j'(t)|$.

Let $\Oh$ denote the double of $\O$ and let $R(z)$ denote
the antiholomorphic involution on $\Oh$ which fixes the boundary of
$\O$.  Let $\Ot=R(\O)$ denote the reflection of $\O$ in $\Oh$ across
the boundary.  We shall frequently use the following fact.  If $g$
and $h$ are meromorphic functions on $\O$ which extend continuously
to the boundary such that
$$g(z)=\overline{h(z)}\qquad\text{for }z\in b\O,$$
then $g$ extends to the double of $\O$ as a meromorphic function.
Indeed, the function $\overline{h(R(z))}$ gives the holomorphic
extension of $g$ to $\Ot$.  For example, since a proper holomorphic
map $f$ from $\O$ to the unit disc extends smoothly up to the boundary
and has modulus one there, it follows that
$$f(z)=1/\overline{f(z)}\qquad\text{for }z\in b\O,$$
and, hence, that $f$ extends to be meromorphic on the double of $\O$.

We now take a moment to recite some standard facts that we shall
assume the reader knows.  Let $A^\infty(\O)$ denote the space of
holomorphic functions on $\O$ that are in $C^\infty(\Obar)$.
Let $L^2(b\O)$ denote the space of complex valued functions on $b\O$
that are square integrable with respect to arc
length measure $ds$.
The Hardy space of functions in $L^2(b\O)$ that are the
$L^2$ boundary values of holomorphic functions on $\O$ shall be written
$H^2(b\O)$.  This space is equal to the closure in $L^2(b\O)$ of
$A^\infty(\O)$ (see \cite{4} for a proof of this elementary fact).

Let $S(z,a)$ denote the classical Szeg\H o kernel associated to the
classical Szeg\H o projection $P$, which is the orthogonal projection of
$L^2(b\O)$ onto $H^2(b\O)$.  For each fixed point $a\in\O$, $S(z,a)$
extends to the boundary as a function of $z$ to be a function in
$A^\infty(\O)$.  Furthermore, $S(z,a)$ has exactly $(n-1)$ zeroes in
$\O$ (counting multiplicities) and does not vanish at any points $z$
in the boundary of $\O$.

The classical Garabedian kernel $L(z,a)$ is a kernel related to the Szeg\H o
kernel via the identity
$$\frac{1}{i} L(z,a)T(z)=S(a,z)\qquad\text{for $z\in b\O$ and $a\in\O$.}
\tag3.1$$
For fixed $a\in\O$, the kernel $L(z,a)$ is a holomorphic function of $z$
on $\O-\{a\}$ with a simple pole at $a$ with residue $1/(2\pi)$.
Furthermore, as a function of $z$, $L(z,a)$ extends to the boundary
and is in the space $C^\infty(\Obar-\{a\})$.  Also, $L(z,a)$ is non-zero
for all $(z,a)$ in $\Obar\times\O$ with $z\ne a$.

The kernel $S(z,w)$ is holomorphic in $z$ and antiholomorphic in $w$
on $\O\times\O$, and $L(z,w)$ is holomorphic in both variables for
$z,w\in\O$, $z\ne w$.  We note here that $S(z,z)$ is real and
positive for each $z\in\O$, and that $S(z,w)=\overline{S(w,z)}$ and
$L(z,w)=-L(w,z)$.  Also, the Szeg\H o kernel reproduces holomorphic functions
in the sense that
$$h(a)=\langle h, S(\cdot,a)\rangle$$
for all $h\in H^2(b\O)$ and $a\in\O$, where the inner product is taken
in $L^2(b\O)$.

Given a point $a\in\O$, the Ahlfors map $f_a$ associated to the pair ($\O,a)$
is related to the Szeg\H o kernel and Garabedian kernel via
$$f_a(z)=\frac{S(z,a)}{L(z,a)}.\tag3.2$$
Note that $f_a'(a)=2\pi S(a,a)\ne 0$.  Because $f_a$ is $n$-to-one, $f_a$
has $n$ zeroes.  The simple pole of $L(z,a)$ at $a$ accounts for the simple
zero of $f_a$ at $a$.   The other $n-1$ zeroes of $f_a$ are given by
$(n-1)$ zeroes of $S(z,a)$ in $\O-\{a\}$.  Let $a_1,a_2,\dots,a_{n-1}$ denote
these $n-1$ zeroes (counted with multiplicity).  It was proved in \cite{5} (see
also \cite{4, page~105}) that, if $a$ is close to one of the boundary
curves, the zeroes $a_1,\dots,a_{n-1}$ become distinct simple zeroes.  It
follows from this result that, for all but at most finitely many points
$a\in\O$, $S(z,a)$ has $n-1$ distinct simple zeroes in $\O$ as a function
of $z$.

The Bergman kernel and the kernel $\Lambda(z,w)$ defined in \S2
satisfy an identity analogous to (3.1):
$$\Lambda(w,z)T(z)=-K(w,z)\overline{T(z)}\qquad\text{for $w\in\O$ and
$z\in b\O$}\tag3.3$$
(see \cite{4, page~135}).  We remark that it follows from well known
properties of the Green's function that
$\Lambda(z,w)$ is holomorphic in $z$ and $w$  and is in
$C^\infty(\Obar\times\Obar-\{(z,z):z\in\Obar\})$.
If $a\in\O$, then $\Lambda(z,a)$ has a double pole at
$z=a$ as a function of $z$ and $\Lambda(z,a)=\Lambda(a,z)$ (see
\cite{4, page~134}).  If $\O$ has real analytic boundary, then
the kernels $K(z,w)$, $\Lambda(z,w)$, $S(z,w)$, and $L(z,w)$,
extend meromorphically to $\Obar\times\Obar$ (see
\cite{4, page~103, 132--136}).
The derivative of the Green's function also satisfies an identity similar
to (3.3):
$$\frac{\dee G}{\dee z}(z,w)T(z)=-
\overline{\frac{\dee G}{\dee z}(z,w)T(z)},\tag3.4$$
for $w\in\O$ and $z\in b\O$ (see \cite{4, page~134}).
We shall also need the identity,
$$S(z,a_1)S(z,a_1)T(z)=-\overline{L(z,a_1)L(z,a_2)T(z)}\tag3.5$$
for $a_1$ and $a_2$ in $\O$ and $z\in b\O$, which follows from (3.1).

The transformation formulas for the Szeg\H o and Garabedian kernels
given in \S2 and the biholomorphic map to a domain with real analytic
boundary allow us to certify that the Ahlfors map is given by formula
(3.2) even in case the domain under study is merely a finitely
connected domain such that no boundary component is a point.

\subhead 4. Proofs of Theorems~2.2 and 2.2a\endsubhead
We first assume that $\O$ is $n$-connected with $n>1$.
Suppose that the Szeg\H o kernel $S(z,w)$ associated to $\O$ is a
rational combination of only two holomorphic functions $G_1$ and $G_2$
on $\O$. To be precise, suppose that $S(z,w)$ is equal to a rational
combination of $G_1(z)$, $G_2(z)$, $\overline{G_1(w)}$, and
$\overline{G_2(w)}$.  It is shown in \cite{8} that
the Bergman kernel $K(z,w)$ associated to $\O$ and all the other
functions mentioned in Theorem~2.2 are rational combinations of
functions of $z$ of the form $S(z,A)$ for three fixed points $A$
in $\O$.  Hence $K(z,w)$ is also a rational combination of the two
holomorphic functions $G_1$ and $G_2$.

Let $f_1$ and $f_2$ denote two Ahlfors maps associated to $\O$ that
generate the field of meromorphic functions on the double of $\O$ (see
\cite{9}).  Let $\C(f_1,f_2)$ denote the field of functions generated
by the two Ahlfors maps (which can be identified with the field of
meromorphic functions on the double in the obvious manner).

It is shown in \cite{7} that any proper holomorphic map from $\O$
onto the unit disc can be expressed as a rational combination
of finitely many functions of $z$ of the form $K(z,w)$ and
$(\dee/\dee\bar w)^mK(z,w)$ for $w$ in $f^{-1}(0)$ (see also
\cite{{11,12}}).  It follows that $f_1$ and $f_2$ are in the field 
$\C(G_1,G_2)$ of functions on $\O$ generated by $G_1$ and $G_2$.
Hence, $\C(f_1,f_2)\subset\C(G_1,G_2)$.

We now claim that the functions $G_1$ and $G_2$ must be algebraically
dependent, i.e., that there exists a complex polynomial $P(z,w)$ of
two complex variables such that $P(G_1(z),G_2(z))\equiv0$.  Indeed
if there did not exist such a polynomial, then we could apply the
Zariski-Castelnuovo Theorem as follows to deduce that $f_1$ and $f_2$
must be algebraically independent, which cannot be the case for
two functions that extend meromorphically to the double of $\O$.
(The Zariski-Castelnuovo Theorem states that an intermediate field
between $\C(x,y)$ and $\C$ must be of the form $\C(u,v)$ or $\C(u)$
where $u$ and $v$ are algebraically independent.)
Since
$$\C(G_1,G_2)\supset\C(f_1,f_2)\supset\C,$$
the algebraic independence of $G_1$ and $G_2$ would imply that
$\C(f_1,f_2)=\C(u,v)$ where $u$ and $v$ are algebraically independent
elements of the field, or $\C(f_1,f_2)=\C(u)$ where $u$ is a single
element of the field.  The first case is clearly impossible because
any two elements of the field of meromorphic functions on the double
are algebraically dependent (see Farkas and Kra \cite{15, p.~248}).
The second case implies that the double of $\O$ is the Riemann sphere,
which only happens if $\O$ is simply connected.

Suppose that $f$ is a proper holomorphic mapping of $\O$ onto the unit
disc.  It is proved in \cite{7} that both $f(z)$ and $f'(z)$ can be expressed
as a rational combinations of finitely many functions of $z$ of the form
$K(z,w)$ and $(\dee/\dee\bar w)^mK(z,w)$ for values of $w$ in the
finite set $f^{-1}(0)$
(see also \cite{{11,12}}).  Hence both $f$ and $f'$ are rational
combinations of $G_1$ and $G_2$.  Since $G_1$ and $G_2$ satisfy a polynomial
identity of the form
$$P(G_1(z),G_2(z))\equiv0$$
on $\O$, we may state that $G_2$ is an algebraic function of $G_1$.
Hence, both $f$ and $f'$ are algebraic functions of $G_1$, and hence
$f'$ must be an algebraic function of $f$, and it follows that
$f$ and $f'$ must be algebraically dependent.

We now consider the simply connected case.
Suppose that the Szeg\H o kernel associated to $\O$ is a
rational combination of only one holomorphic function $G$.  Since
the Bergman kernel is a constant times the square of the Szeg\H o
kernel, it follows that it is also a rational combination of $G$.
Let $f$ denote a Riemann map associated to $\O$ mapping $\O$ one-to-one
onto the unit disc.  Note that $f$ extends to the double of $\O$ and that
it generates the field of meromorphic functions on the double of $\O$ (see
\cite{15}).

The same proof given in \cite{7} for proper maps shows that $f(z)$ and
$f'(z)$ are both rational combinations of functions of $z$ of the form
$K(z,w)$ and $(\dee/\dee\bar w)K(z,w)$ and $w=f^{-1}(0)$.
Hence both $f$ and $f'$ are rational combinations of $G$.  Hence $f'$
must be an algebraic function of $f$, and it follows that $f$ and $f'$
must be algebraically dependent.

To finish the proofs of Theorems~2.2 and 2.2a, we now turn to the
business of constructing a special Riemann surface attached to the domain
in case a proper map onto the unit disc and its derivative are algebraically
dependent.

\subhead 5. Construction of a Riemann surface\endsubhead
We assume at first that $\O$ is a finitely connected domain in the
plane with $C^\infty$ smooth boundary.  Suppose also that there exits
a proper holomorphic mapping $f$ from $\O$ onto the unit disc with the
property that $f$ and $f'$ are algebraically dependent, i.e.,
there exists a polynomial $P(z,w)$ such that $P(f'(z),f(z))\equiv0$
on $\O$.  We now construct the Riemann surface mentioned in
Theorems~2.2 and~2.2a.

Let $S_{\bar m}(z,a)$ denote $(\dee^m/\dee\bar a^m)S(z,a)$ and
let $L_m(z,a)$ denote $(\dee^m/\dee a^m)L(z,a)$.  Notice that we may
differentiate the conjugate of (3.1) $m$-times with respect to $\bar a$
to obtain
$$S_{\bar m}(z,a)=
i \overline{L_m(z,a)T(z)}$$
for $z\in b\O$ and $a\in\O$.  When the square of this identity is
combined with the identity
$$\frac{f'(z)}{f(z)}T(z)=-\overline{f'(z)T(z)/f(z)}
\quad\text{for }z\in b\O$$
(which is obtained by differentiating $\log|f(z(t))|=0$ with
respect to $t$ when $z(t)$ is a parameterization of a boundary
curve), we see that 
$$\frac{S_{\bar m}(z,a)^2f(z)}{f'(z)}$$
is equal to the conjugate of
$$\frac{L_m(z,a)^2f(z)}{f'(z)}$$
for $z\in b\O$ and any point $a$ in $\O$.
This identity reveals that
$S_{\bar m}(z,a)^2f(z)/f'(z)$ extends to the double of $\O$ as a meromorphic
function of $z$.  Since $f(z)$ extends to the double of $\O$, it follows that
$S_{\bar m}(z,a)^2/f'(z)$ extends to the double of $\O$ as a meromorphic
function.  The algebraic dependence of $f$ and $f'$ means that they
satisfy a polynomial equation $P(f'(z),f(z))\equiv0$ on $\O$.  Since
$f(z)$ extends to the double of $\O$ as a meromorphic function, this
polynomial identity reveals how to extend $f'(z)$ to the double of $\O$
as a finitely valued function with at most finitely many algebraic
singularities.  Furthermore, since $S_{\bar m}(z,a)^2/f'(z)$ extends to
the double of $\O$ as a meromorphic function, we may state that
$S_{\bar m}(z,a)$ extends to the double of $\O$  as a finitely valued
function of $z$ with at most finitely many algebraic singularities for
each fixed point $a$ in $\O$.

Theorem~2.6 states that the Szeg\H o kernel is given by a rational
combination of $f$  and finitely many functions
of the form $S_{\bar m}(z,a_j)$.  Let us call the functions in the
list of finitely many functions of the form $S_{\bar m}(z,a_j)$
just mentioned {\it core functions}.
We may view the core functions (from the viewpoint
of Weierstrass) as being finitely valued multivalued functions that can
be analytically continued to the double of $\O$.  There is a finite set
of points $E$ in double of $\O$ at which one or more of the function
elements associated to the core functions has an algebraic singularity.
Choose a point $A_0$ in $\O-E$ to act as a base point.
We construct $\Cal R$ by performing analytic continuation
of each of the core functions {\it simultaneously}, starting at $A_0$
and moving all around $\Oh$, paying special attention to
the points in $E$.  Away from $E$, the lifting of germs along curves
to a Riemann surface over the double of $\O$ is routine and obvious.
When we analytically continue up to a point $p$ in $E$, it may happen that
none of the germs of the function elements of the core functions become
singular at $p$.  In this case, we lift and analytically continue through $p$
without incident.  If, on the other hand, at least one of the elements
is singular at $p$, then we construct a local coordinate system at
the point $\tilde p$ above $p$ as follows.  Consider the function
elements of the core functions that are obtained as we analytically
continue them up to $p$ along a curve.  Each of these elements can be
viewed as a function element of a Puiseux expansion at $p$ in a local
coordinate $\zeta$ where $p$ corresponds to the origin.  Hence, for each
core function $C(z)$ there is a positive integer $\lambda$ such that the
substitution $z=p+(\zeta)^{\lambda}$ makes $C(\zeta)$ analytic and continuable
in $\zeta$ through $\zeta=0$.  (Note that the number $\lambda$ is equal to
one if $C(z)$ does not have a singularity at $p$.)
Let $m$ be equal to the least common multiple of all the numbers
$\lambda$ associated to the core functions.
We can now define a local uniformizing variable that is
suitable for each of the function elements in the obvious manner:
$z=p+(\zeta)^{m}$.  This coordinate function allows us to lift all the
core functions so as to be defined and single valued on a disc centered
at $\zeta=0$ and we use it to define a local chart near $\tilde p$.

It is clear that we can identify $\O$ as a subdomain of the Riemann surface
$\Cal R$ by virtue of the fact that the core functions all have {\it
preferred\/} germs in $\O$ given by the values they have via the definition
of the kernel functions on the domain $\O$.  A point $\tilde{p}\in{\Cal R}$
over $p\in\O$ shall be identified as the point $p$ in $\O\subset{\Cal R}$ if
the germs of all the core functions are equal to their preferred germs at $p$.
Note that there will be other sheets of $\Cal R$ above $\O$ if it
happens that one or more of the core functions continue back to $p\in\O$
and are not equal to their preferred germs over $\O$.
The Riemann surface $\Cal R$ is clearly compact because of the finite
valuedness of the core functions on the compact $\Oh$.  Note that any
meromorphic function on $\Oh$ can be defined as a meromorphic
function on $\Cal R$.  Hence, $f(z)$ can be extended to $\Cal R$ as
a meromorphic function.  Since the core functions and $f$ extend to
$\Cal R$, it follows that the Szeg\H o kernel extends to be
meromorphic on $\Cal R\times\Cal R$.  Since the Bergman kernel is a rational
combination of three functions of $z$ of the form $S(z,A)$, it follows
that the Bergman kernel extends to be meromorphic on $\Cal R\times\Cal R$.
Similarly, so do the functions $F_j'$.

To see that the complement of $\O$ in $\Cal R$ is connected, we can
follow exactly the same pole counting procedure given in \cite{10} where
a similar Riemann surface is attached to a domain with algebraic kernel
functions.  This procedure yields that each boundary
component of $\O$ is attached to the Riemann surface exactly once and that
the complement of $\O$ in $\Cal R$ is connected.

Let $G_1$ and $G_2$ denote a primitive pair for $\Cal R$.  Theorem~2.6
states that the Szeg\H o kernel is a rational combination of $f(z)$ and
the core functions.  Since $f(z)$ and the core functions extend to be
meromorphic functions on $\Cal R$, we may now state that the Szeg\H o
kernel is a rational combination of $G_1$ and $G_2$.  Since the Bergman
kernel associated to $\O$ and all the other functions mentioned in
Theorem~2.2 are rational combinations of functions of $z$ of the form
$S(z,A)$ for three fixed points $A$ in $\O$, they too are rational
combinations of $G_1$ and $G_2$.  The proof of Theorem~2.2 is complete
in case $\O$ has smooth boundary.

If $\O$ does not have smooth boundary, we shall use the identity
$$P(f(z),f'(z))\equiv0,\tag5.1$$
which is assumed to hold for a proper holomorphic mapping $f(z)$ of
$\O$ onto the unit disc, to show that the boundary of $\O$ must be given
by piecewise $C^\infty$ smooth real analytic curves.  Let $F(w)$ denote a
local inverse defined on a small open subset of $D_1(0)$ to the proper
holomorphic map $f(z)$.  If we replace $z$ by $F(w)$ in (5.1) we obtain
$P(w,1/F'(w))\equiv0$, and this shows that $F'(w)$ is an algebraic
function.  Hence $F$ can be analytically continued past the boundary
of the unit disc except at possibly finitely many points where $F'$ has an
algebraic singularity.  Since the continuation of $f$ maps the boundary of
the unit disc into the boundary of $\O$, it is easy to deduce that the
boundary of $\O$ is given by piecewise real analytic curves.  (Only
minor complications are introduced if the point at infinity is in the
boundary of $\O$ and these are bypassed by standard arguments.)  Now
the same construction of $\Cal R$ we used above can be carried out and
the pole counting argument of \cite{10} given in the case of piecewise
real analytic boundary applies to yield that the complement of $\O$ in
$\Cal R$ is connected.  The rest of the proof is the same.

\subhead 6. The proofs of Theorems 2.3 and 2.5\endsubhead
We shall see momentarily that we will be able to reduce our problem
to the case where $\O$ has $C^\infty$ smooth real analytic boundary.
The following lemmas will allow us to see that the two classes~$\Cal A$
and~$\Cal B$ are natural and that they are the largest classes of
functions that can appear in the statements of Theorems~2.3 and~2.5.

\proclaim{Lemma 6.1}
Suppose that $\O$ is a finitely connected domain in the plane with
$C^\infty$ smooth real analytic boundary.
On such a domain, the class~$\Cal B$ is equal to
the set of meromorphic functions $G$ on $\Obar$ that satisfy an
identity of the form
$$G(z)=\overline{H(z)T(z)}\tag6.1$$
for $z\in b\O$, where $H(z)$ is another meromorphic function on $\Obar$.
\endproclaim

When we deal with the class~$\Cal A$ we shall need the following lemma.

\proclaim{Lemma 6.2}
Suppose that $\O$ is a finitely connected domain in the plane with
$C^\infty$ smooth real analytic boundary.
On such a domain, the class~$\Cal A$ is equal to
the set of meromorphic functions $G$ on $\Obar$ that satisfy an
identity of the form
$$G(z)T(z)=\overline{H(z)T(z)}\tag6.2$$
for $z\in b\O$, where $H(z)$ is another meromorphic function on $\Obar$.
\endproclaim

These two lemmas have the important consequence that if $g_1/g_2$ is a
quotient of two functions in the class $\Cal A$ (or two functions in
the class $\Cal B$) where $g_2\not\equiv0$, then $g_1/g_2$ extends to
the double of $\O$ as a meromorphic function.  Indeed, this fact
follows directly from the identity of the form
$g_1(z)/g_2(z)=\overline{h_1(z)/h_2(z)}$ on $b\O$ that would
be satisfied by the quotient of two such functions.

\demo{Proof of Lemma 6.1}
Notice that (3.1) shows that $L(z,a)$ and $S(z,a)$
satisfy the similar identities:
$$\gather
S(z,a)=\overline{-iL(z,a)T(z)}\\
L(z,a)=\overline{-iS(z,a)T(z)}.
\endgather$$
Furthermore, these identities can be differentiated with respect to
$a$ or $\bar a$.
$$\gather
S_{\bar m}(z,a)=\overline{-iL_m(z,a)T(z)}\\
L_m(z,a)=\overline{-iS_{\bar m}(z,a)T(z)}.
\endgather$$
Hence, functions in the class $\Cal B$ satisfy the
condition of formula (6.1).

Now suppose that $G$ and $H$ are meromorphic functions on $\Obar$ that satisfy
$$G(z)=\overline{H(z)T(z)}$$
for $z\in b\O$.  When $a$ is in $\O$, we note that $L(z,a)$ has a
single simple pole at $a$ as a function of $z$ and $S(z,a)$ has no
singularities in $\Obar$, and $L_m(z,a)$ has a single pole of order~$m+1$
at $a$ and $S_{\bar m}(z,a)$ has no singularities in $\Obar$.  Hence, we may
subtract linear combinations of the identities in the paragraph above from
(6.1) designed to remove the singularities of $G$ and $H$ in $\O$ from both
sides of the equation.  In the end, we obtain an identity of the form
$$G_1(z)=\overline{H_1(z)T(z)}$$
where $G_1$ and $H_1$ are holomorphic on $\O$ with possibly finitely many
poles on $b\O$.   When $a$ is in the boundary, the function $S(z,a)$
has a single simple pole at $a$ in $\Obar$ and $S_{\bar m}(z,a)$ has a single
pole of order~$m+1$ at $a$ in $\Obar$.  Hence, we can subtract linear
combinations of $S(z,a)$ and $S_{\bar m}(z,a)$ from $G_1$ to eliminate the
poles of $G_1$ in $\Obar$ and subtract the corresponding linear
combinations of $-iL(z,a)$ and $-iL_m(z,a)$ from $H_1$ to obtain an
identity of the form
$$g(z)=\overline{h(z)T(z)}$$
where $g$ is holomorphic on $\Obar$ and $h$ is holomorphic on $\O$ with
possibly finitely many poles on $b\O$.  However, the identity
$g(z)=\overline{h(z)T(z)}$ on $b\O$ shows that $h$ does not tend to
infinity at any point on the boundary, so $h$ actually has no poles on
$\Obar$.  Hence, both $g$ and $h$ are holomorphic on $\Obar$.  Now we can
conclude from the fact that $\overline{h(z)T(z)}$ is orthogonal
to holomorphic functions in $L^2(b\O$) that $g$ and $h$ must both be
zero.  The proof of our claim is complete.
\enddemo

\demo{Proof of Lemma 6.2}
The three identities
$$\gather
F_j'(z)T(z)=-\overline{F_j'(z)T(z)}\\
G_z(z,a)T(z)=-\overline{G_z(z,a)T(z)}\\
S(z,a_1)S(z,a_2)T(z)=-\overline{L(z,a_1)L(z,a_2)T(z)}
\endgather$$
for $z\in b\O$ (see \cite{4, pages~80,~134}, (3.3), (3.4), and (3.5)) reveal
that functions in the class $\Cal A$ satisfy the condition of formula (6.2).

Now suppose that $G$ and $H$ are meromorphic functions on $\Obar$ that satisfy
$$G(z)T(z)=\overline{H(z)T(z)}$$
for $z\in b\O$.  Note that the function $G_z(z,a)$ has a simple pole at $a$
even when $a$ is in $b\O$.  Differentiate the identity
$G_z(z,a)T(z)= -\overline{G_z(z,a)T(z)}$ ($z\in b\O$)
$m$~times with respect to $a$.  The function
$(\dee^m/\dee a^m)G_z(z,a)$ on the left hand side of the formula is
a meromorphic function of $z$ on $\Obar$ with a single pole of order $m+1$
at $z=a$.  This is true even if $a\in b\O$.  When $a\in\O$, the function
$(\dee^m/\dee\bar a^m)G_z(z,a)$ on the right hand side of the formula is
a holomorphic function of $z$ on $\Obar$ with no
singularity at $z=a$.  If we differentiate the same formula $m$~times
with respect to $\bar a$, then the function
$(\dee^m/\dee\bar a^m)G_z(z,a)$ on the left hand side of the formula is
a holomorphic function of $z$ on $\Obar$ with no
singularity at $z=a$ when $a\in\O$ and the function
$(\dee^m/\dee a^m)G_z(z,a)$ on the right hand side of the formula is
a meromorphic function of $z$ on $\Obar$ with a single pole of order $m+1$
at $z=a$.  Hence, it is possible to subtract linear
combinations of functions in the class~$\Cal A$ from the identity
$GT=\overline{HT}$ to remove all the poles of $G(z)$ in $\Obar$ and
all the poles of $H$ of order two or more in $\O$ to obtain an
identity of the form
$$g(z)T(z)=\overline{h(z)T(z)}$$
where $g$ is holomorphic on $\Obar$ and $h$ is meromorphic on $\Obar$
with only finitely many simple poles in $\O$ and finitely many other
poles on $b\O$.  However, none of the poles of $h$ can actually occur on
the boundary because the identity $gT=\overline{hT}$ shows that $h$
cannot blow up there.  Let $\{b_j\}_{j=1}^N$ denote the simple poles
of $h$ in $\O$ and choose a point $a$ in $\O$ distinct from these
points.  The functions $S(z,b_j)S(z,a)$ are holomorphic on $\Obar$
and satisfy the identity
$$S(z,b_j)S(z,a)T(z)=-\overline{L(z,b_j)L(z,a)T(z)}$$
for $z\in b\O$.  Since $L(z,b_j)L(z,a)T(z))$ has a simple pole at
$b_j$, we may subtract a linear combination of these identities from
$g(z)T(z)=\overline{h(z)T(z)}$ to eliminate the poles of $h$ at the
points $b_j$.  We obtain
$$g_1(z)T(z)=\overline{h_1(z)T(z)}$$
where $g_1$ is holomorphic on $\Obar$ and $h_1$ is meromorphic on
$\Obar$ with possibly a single simple pole at $a$.  However, if
$g_1(z)T(z)=\overline{h_1(z)T(z)}$ is integrated around the boundary
with respect to arc length measure, the left hand side is zero by
Cauchy's Theorem and the right hand side is equal to $2\pi i$ times
the residue of $h_1$ at $a$.  So, in fact, $h_1$ has no poles in
$\Obar$.  It is shown in \cite{4, page~80} that holomorphic
functions $g_1$ and $h_1$ that satisfy
$g_1(z)T(z)=\overline{h_1(z)T(z)}$ on the boundary must be linear
combinations of the functions $F_j'$, and this shows that $G$ is in
the class~$\Cal A$ and the proof of our claim is complete.
\enddemo

Recall that the classes $\Cal A$ and $\Cal B$ were defined
differently for domains with smooth boundary than for domains
without.  We now let $\Cal A$ and $\Cal B$ denote
the classes defined by restricting the points $a$ in the
definitions of the classes $\Cal A$ and $\Cal B$ to be in
$\O$.  If $\O$ has smooth boundary, then we let
${\Cal A}^+$ and ${\Cal B}^+$ denote the classes obtained
by letting the points $a$ in the definitions also fall on the
boundary.

\proclaim{Lemma 6.3}
Suppose that $\Phi:\O\to\O_a$ is a biholomorphic mapping between a domain $\O$
and a finitely connected domain $\O_a$ with $C^\infty$ smooth real analytic
boundary such that no boundary component is a point.  The
transformation
$$h\mapsto \Phi'(z)h(\Phi(z))$$
is a one-to-one linear map of the
class~$\Cal A$ associated to $\O_a$
onto
the class~$\Cal A$ associated to $\O$.
If $b\O$ is $C^\infty$ smooth,
then this transformation is also a one-to-one
linear map of the class ${\Cal A}^+$ associated to $\O_a$
onto
the class~${\Cal A}^+$ associated to $\O$.
The transformation
$$h\mapsto \sqrt{\Phi'(z)}h(\Phi(z))$$
is a one-to-one linear map of the
class~$\Cal B$ associated to $\O_a$
onto
the class~$\Cal B$ associated to $\O$.
If $b\O$ is $C^\infty$ smooth,
then this transformation is also a one-to-one
linear map of the class ${\Cal B}^+$ associated to $\O_a$
onto
the class~${\Cal B}^+$ associated to $\O$.
\endproclaim

The proof of Lemma~6.3 is straightforward and uses only the well
known transformation properties of the functions that generate the two
classes.  We omit the proof.

\demo{Proof of Theorem 2.5}
Assume for the moment that $\O$ has smooth real analytic boundary.
Fix a point $a$ in $\O$ so that the zeroes $a_1,\dots,a_{n-1}$ of
$S(z,a)$ are distinct simple zeroes.  I proved in \cite{4, Theorem~3.1}
that the Szeg\H o kernel can be expressed via
$$S(z,w)=\frac{1}{1-f_a(z)\overline{f_a(w)}}\left(c_0
S(z,a)\overline{S(w,a)}+
\sum_{i,j=1}^{n-1} c_{ij}S(z,a_i)\,\overline{S(w,a_j)}\right)\tag6.3$$
where $f_a(z)$ denotes the Ahlfors map associated to $(\O,a)$,
$c_0=1/S(a,a)$, and the coefficients $c_{ij}$ are given as the
coefficients of the inverse matrix to the invertible matrix
$\left[S(a_j,a_k)\right]$.  This shows that
$$S(z,w)(1-f_a(z)\overline{f_a(w)})$$
is a linear combination of functions of the form $g(z)\overline{h(w)}$
where $g$ and $h$ are in the class~$\Cal B$.  Let $B(z)$ be any non-zero
function in the class~$\Cal B$.  The remark after the statement of Lemma~6.2
asserts that the quotient of any two functions in the class~$\Cal B$ extends
to the double of $\O$ as a meromorphic function.  It follows that
$$\frac{S(z,w)(1-f_a(z)\overline{f_a(w)})}{B(z)\overline{B(w)}}
\tag6.4$$
is equal to a linear combination of functions of the form $g(z)/B(z)$ times
the conjugate of $h(w)/B(w)$, and because these quotients are mermorphic on
the double, they can be expressed as rational functions of the
primitive pair.  The Ahlfors map itself extends to the double as a
meromorphic function, and hence the claim about the Szeg\H o kernel
is proved.

In case $\O$ does not have smooth real analytic boundary, we use the
well known fact that there exists a biholomorphic map $\Phi$ from $\O$ onto
a bounded domain $\O_a$ with smooth real analytic boundary.  The
transformation formulas for the Szeg\H o kernel, the Ahlfors maps,
and the functions in the class~$\Cal B$ reveal that formula (6.4) transforms
by simple composition with $\Phi$, i.e., the terms involving $\Phi'$
cancel.  Similarly, the linear combination of terms of the form
$g(z)/B(z)$ times the conjugate of $h(w)/B(w)$ where $g$ and $h$
are in the class~$\Cal B$ enjoy the same property, and these extend to the
double as meromorphic functions.  The proof is complete.
\enddemo

\demo{Proof of Theorem~2.3}
Assume for the moment that $\O$ has smooth real analytic boundary.
The Bergman kernel $K(z,w)$ is related to the Szeg\H o kernel via the
identity
$$K(z,w)=4\pi S(z,w)^2+\sum_{i,j=1}^{n-1}
A_{ij}F_i'(z)\overline{F_j'(w)},\tag6.5$$
Notice that formula (6.3) and Lemma~6.2 yield that $S(z,w)^2$ is a linear
combination of functions in the class~$\Cal A$ divided by
$(1-f_a(z)\overline{f_a(w)})^2$.  Hence, if we divide (6.5) by
$A(z)\overline{A(w)}$ where $A(z)$
is a non-zero function in the class~$\Cal A$, then we see that
$$\frac{K(z,w)}{A(z)\overline{A(w)}}$$
is a linear combination of functions of the form $g(z)/A(z)$ times
conjugates of $h(w)/A(w)$ times either constants or
$(1-f_a(z)\overline{f_a(w)})^{-2}$, where $g$ and $h$ are in the
class~$\Cal A$.  Since all these functions of $z$ and $w$ extend to the
double of $\O$ as meromorphic functions, Theorem~2.3 follows.  Now
the same argument given in the last paragraph of the proof of
Theorem~2.5 shows that we can drop the hypothesis that the boundary
of $\O$ be smooth and real analytic.  The proof is complete.
\enddemo

\subhead 7. The Szeg\H o kernel with weights \endsubhead
Suppose that $\O$ is a bounded $n$-connected domain in the
plane with $C^\infty$ smooth boundary, i.e., a domain whose boundary
$b\O$ is given by finitely many non-intersecting $C^\infty$ simple closed
curves.  Recall that $L^2(b\O)$ denotes the space of complex valued functions
on $b\O$ that are square integrable with respect to arc length measure $ds$.
Given a positive real valued $C^\infty$ function $\phi$ on the boundary
of $\O$, let $L^2_\phi(b\O)$ denote the space of complex valued functions
on $b\O$ that are square integrable with respect to $\phi(s)ds$.
The Hardy space of functions in $L^2(b\O)$ that are the
$L^2$ boundary values of holomorphic functions on $\O$ shall be written
$H^2(b\O)$.  This space is equal to the closure in $L^2(b\O)$ of
$A^\infty(\O)$ (see \cite{4}) and can be identified in a natural way
with the subspace of $L^2_\phi(b\O)$ equal to the closure of $A^\infty(\O)$ in
that space.  Thus, we need not define $H^2_\phi(b\O)$ separately.

The inner products associated to $L^2(b\O)$ and $L^2_\phi(b\O)$
shall be written
$$
\langle u,v\rangle=\int_{b\O}u\ \bar v\ ds,
\quad\text{ and }\quad
\langle u,v\rangle_{\phi}=\int_{b\O}u(s)\ \overline{v(s)}\ \phi(s)ds,$$
respectively.
We let $S(z,a)$ denote the classical Szeg\H o kernel associated to the
classical Szeg\H o projection $P$, which is the orthogonal projection of
$L^2(b\O)$ onto $H^2(b\O)$, and we let $\sigma(z,w)$ denote the kernel
associated to the orthogonal projection $P_\phi$ of the weighted space
$L^2_\phi(b\O)$ onto $H^2(b\O)$.  The arguments presented in \cite{4}
showing that $P$ maps $C^\infty(b\O)$ into itself can easily be modified
to show that $P_\phi$ also maps $C^\infty(b\O)$ into itself.  Also, if
the boundary curves of $\O$ are $C^\infty$ smooth {\it real analytic\/}
curves and $\phi$ is real analytic on $b\O$, then $P_\phi$ maps real analytic
functions on $b\O$ into the space of holomorphic functions on $\O$ that
extend to be holomorphic on a neighborhood of $\Obar$ (see \cite{4, page~41}).

Any function $v$ in the subspace $H^2(b\O)^\perp$ of $L^2(b\O)$ (which is
the orthogonal complement of $H^2(b\O)$ in $L^2(b\O)$) can be written
$$v= \overline{HT},$$
for a unique $H$ in $H^2(b\O)$ (see \cite{4, p.~13}).  Consequently,
every function $u$ in $L^2(b\O)$ has an orthogonal decomposition of
the form
$$u=h+\overline{HT},$$
where $h\in H^2(b\O)$ and $H\in H^2(b\O)$.  We can easily deduce from
this fact that every function $u$ in
$L^2_\phi(b\O)$ has an orthogonal decomposition of the form
$$u=h+\phi^{-1}\overline{HT},$$
where $h\in H^2(b\O)$ and $H\in H^2(b\O)$.  Indeed, let $h=P_\phi u$.
Then $u-h$ is orthogonal to $H^2(b\O)$ with respect to the weighted
inner product.  Hence $(u-h)\phi$ is orthogonal to $H^2(b\O)$ with
respect to the standard inner product on $L^2(b\O)$.  Hence
$(u-h)\phi=\overline{HT}$ for a unique $H\in H^2(b\O)$ and our claim
is proved.  This orthogonal decomposition allows us to define a
{\it weighted Garabedian kernel\/} as follows.  The Cauchy integral
formula reveals that the {\it weighted Cauchy kernel\/} $C_a(z)$, which
is given as the conjugate of
$$\frac{1}{2\pi i}\frac{\phi(z)^{-1}T(z)}{z-a},$$
reproduces holomorphic functions with respect to the weighted inner
product in the sense that
$$h(a)=\langle h,C_a\rangle_\phi.$$
Hence, it follows that
$$\sigma(z,a)=P_\phi C_a$$
and the orthogonal decomposition for $C_a$ is
$$C_a(z)=\sigma(z,a) +\phi^{-1}\overline{H_aT},$$
where $H_a$ is in $H^2(b\O)$.
By analogy with the definition of the classical Garabedian kernel
(see \cite{4, p.~24}), we define the weighted Garabedian kernel
$\lambda(z,a)$ to be given by
$$\lambda(z,a)=\frac{1}{2\pi }\frac{1}{z-a}-iH_a(z).$$
Notice that $\sigma(z,a)$ and $\lambda(z,a)$ satisfy the identity
$$\overline{\sigma(z,a)}=\frac{1}{i\phi(z)}\lambda(z,a)T(z)
\tag7.1$$
for $z\in b\O$ and $a\in\O$.  Since $\sigma(z,a)=P_\phi C_a$, it
follows that $\sigma(z,a)$ is in $A^\infty(\O)$ as a function of $z$
for each fixed $a\in\O$ if $\phi$ is $C^\infty$ smooth.  Furthermore,
for a fixed point $a\in\O$, $\lambda(z,a)$ is holomorphic function of $z$
on $\O-\{a\}$ with a simple pole at $a$ with residue $1/(2\pi)$ and
$\lambda(z,a)$ extends $C^\infty$ smoothly to $b\O$ if $\phi$ is
$C^\infty$ smooth.  In case the boundary of $\O$ is real analytic
and $\phi$ is real analytic on $b\O$, both $\sigma(z,a)$ and
$\lambda(z,a)$ extend holomorphically past the boundary in $z$ for
each fixed $a$ in $\O$.  We record here for future use the
derivative
$$\frac{\dee^n}{\dee a^n}\sigma(a,z)=
\frac{1}{i\phi(z)}\frac{\dee^n}{\dee a^n}\lambda(z,a)T(z)
\tag7.2$$
of (7.1) with respect to $a$.

The weighted Szeg\H o kernel reproduces holomorphic functions with
respect to the weighted inner product in the sense that
$$h(a)=\langle h, \sigma(\cdot,a)\rangle_{\phi}$$
for $h\in H^2(b\O)$.  This last identity may be differentiated with
respect to $a$ to yield that
$$h^{(n)}(a)=\langle h, (\dee^n/\dee\bar a^n)\sigma(\cdot,a)\rangle_{\phi}$$
for $h\in H^2(b\O)$.  The weighted Szeg\H o kernel satisfies
$\sigma(z,a)=\overline{\sigma(a,z)}$.

The conjugate of formula (7.1) is
$$\frac{i}{\phi(z)}\overline{\lambda(z,a)}=\sigma(z,a)T(z)$$
If we multiply this by (7.1) and divide out the $1/\phi$ factors, we obtain
$$\sigma(z,a)\lambda(z,a)T(z)=-
\overline{\lambda(z,a)\sigma(z,a)T(z)}\quad\text{for }z\in b\O$$
and $a\in\O$.  Hence, it follows from Lemma~6.2 that
$\sigma(z,a)\lambda(z,a)$ is in the class $\Cal A$ and hence, even though the
weight function may be quite arbitrary, it is possible to relate these
kernels to standard objects of potential theory.  For example, similar
arguments to those used in \cite{8} can be used to show that the Green's
function can be expressed via
$$\frac{\dee}{\dee z} G(z,w)=
\pi\frac{\sigma(z,w)\lambda(z,w)}{\sigma(w,w)}+
\sum_{j=1}^{n-1}c_j(w)F_j'(z)$$
where the functions $c_j(w)$ can be easily determined as in \S5 of \cite{8}
to be given by a linear combination of $n-1$ explicit harmonic functions and
finitely many other functions that are rational combinations of the
basic holomorphic functions that comprise $\sigma$ and $\lambda$.

\subhead 8. An orthonormal basis for the weighted Hardy space \endsubhead
We suppose, as we did in \S7, that $\O$ is a bounded $n$-connected domain
in the plane with $C^\infty$ smooth boundary and that $\phi$ is a real
valued positive $C^\infty$ weight function on $b\O$.   Suppose that
$f:\O\to D_1(0)$ is a proper holomorphic mapping of $\O$ into the unit disc.
Such a mapping extends $C^\infty$ smoothly to the boundary and is a finite
branched covering map of some finite order $M$.
Let $a_1,\dots,a_N$ denote the zeroes of $f$ in $\O$ and let $M(k)$
denote the multiplicity of the zero of $f$ at $a_k$.  Of course,
$M=\sum_{k=1}^N M(k)$.  Notice that, because $|f(z)|=1$ for $z$ in
$b\O$, it follows that  $f(z)=1/\,\overline{f(z)}$ for $z$ in $b\O$.

Let $\sigma_0(z,w)$ denote $\sigma(z,w)$, let
$\sigma_{\bar n}(z,w)$ denote $(\dee^n/\dee\bar w^n)\sigma(z,w)$,
and let
$$\sigma_{m\bar n}(z,w):=\frac{\dee^{m+n}}{\dee z^m\dee\bar w^n}\sigma(z,w).$$
We shall now prove that the set of functions
$$h_{inp}(z)=\sigma_{\bar n}(z,a_i)f(z)^p,$$
where $1\le i\le N$, $0\le n\le M(i)$, and $p\ge 0$,
forms a basis for the Hardy space $H^2(b\O)$ and that
$$\langle h_{inp},h_{jmq}\rangle_{\phi}=
\cases
0, &\text{if $p\ne q$}, \\
\sigma_{m\bar n}(a_j,a_i),
&\text{if $p=q$}.\endcases\tag8.1$$

First, we must show that the set of functions above spans a dense subset of
$H^2(b\O)$.  Indeed, suppose that $g\in H^2(b\O)$ is orthogonal to the span.
Notice that the reproducing property of the weighted Szeg\H o kernel yields
that
$$0=\langle g,\sigma_{\bar n}(\cdot,a_j)\rangle_{\phi}=g^{(n)}(a_j)$$
for $0\le n\le M(j)$,
and therefore $g$ vanishes at $a_1,\dots,a_N$ to the same order that $f$ does.
Suppose we have shown that $g$ vanishes to order $m$ times the order that $f$
vanishes at each $a_j$, $j=1,\dots,N$.  It follows that $g/f^m$ has removable
singularities at each $a_j$ and so it can be viewed as an element of
$H^2(b\O)$.  We shall now show that $g/f^m$ must vanish to the same order
at each $a_j$ that $f$ does.  Since $1/f(z)=\overline{f(z)}$ when $z\in b\O$,
we may write
$$0=\langle g,\sigma_{\bar n}(\cdot,a_j)f^m\rangle_{\phi}=
\langle g/f^m,\sigma_{\bar n}(\cdot,a_j)\rangle_{\phi}$$
and this last quantity is equal to the $n$-th derivative of $g/f^m$ at $a_j$.
Since this is zero for $0\le n\le M(j)$, we conclude that $g/f^m$ vanishes
to the same order that $f$ does at each $a_j$.  Hence, $g$ vanishes to order
$m+1$ times the order that $f$ vanishes at each $a_j$, $j=1,\dots,N$.  By
induction, we may conclude that $g$ vanishes to infinite order at each $a_j$
and hence, $g\equiv0$.  This proves the density.

We now turn to the proof of (8.1).  We may suppose that $p\ge q$.  The fact
that $\overline{f}=1/f$ on $b\O$ yields that
$$\langle h_{inp},h_{jmq}\rangle_{\phi}=
\langle\sigma_{\bar n}(z,a_i)f(z)^{p-q},\sigma_{\bar m}(z,a_j)\rangle_{\phi}.$$
The reproducing property of the weighted Szeg\H o kernel yields that
this last last inner product is equal to
$$\frac{\dee^m}{\dee z^m}\left[\sigma_{\bar n}(z,a_i)f(z)^{p-q}\right]$$
evaluated at $z=a_j$.
Since the multiplicity of the zero of $f$ at $a_j$ is greater than or equal to
$m$, this quantity is zero if $p>q$.  If $p=q$, then the
$f(z)$ term is not present and the proof of identity (8.1) complete.

It is now easy to see that
the functions $h_{inp}$ are linearly independent.  Indeed, identity
(8.1) reveals that we need only check that, for fixed $p$,  the
functions $h_{inp}$, $i=1,\dots,N$, $n=0,\dots,M(i)$, are linearly
independent, and this is true because a relation of the form
$$\sum_{i=1}^{N}\sum_{n=0}^{M(i)} C_{in} \sigma_{\bar n}(z,a_i)\equiv 0$$
implies, via the reproducing property of the weighted Szeg\H o kernel, that
every function $g$ in the Hardy space satisfies
$$\sum_{i=1}^{N}\sum_{n=0}^{M(i)} \overline{C_{in}} g^{(n)}(a_i)=0,$$
and it is easy to construct polynomials $g$ that violate such a
condition.

We next orthonormalize the sequence $\{h_{inp}\}$ via the Gram-Schmidt
procedure.  Formula (8.1) shows that we need only
orthonormalize the functions $h_{inp}$, $i=1,\dots,N$,
$n=0,\dots,M(i)$ for each fixed $p$.  If we are careful to perform the
Gram-Schmidt procedure in exactly the same order with respect to the
indices $i$ and $n$ at each level $p$, we obtain an orthonormal set
$\{H_{inp}\}$ which is related to our original set via a formula,
$$
H_{inp}(z)=\sum_{j=1}^{N}\sum_{m=0}^{M(j)}
b_{ijnm} h_{jmp},
$$
where, because $|f|=1$ on $b\O$,
{\it the coefficients\/} $b_{ijnm}$ {\it do not depend on\/}
$p$.  This last fact is critical in what follows.

The weighted Szeg\H o kernel can be written in terms of our orthonormal basis
via
$$\sigma(z,w)=\sum_{p=0}^\infty
\sum_{i=1}^{N}
\sum_{n=1}^{M(i)}
 H_{inp}(z)\,\overline{H_{inp}(w)}.$$
The geometric sum
$$\sum_{p=0}^\infty f(z)^p\,\overline{f(w)^p}=\frac{1}{1-f(z)\overline{f(w)}}$$
can be factored from the expression for $\sigma(z,w)$ to yield a formula like
the one in the following theorem.  (Note that by taking the weight
function to be identically one, we obtain a proof of Theorem~2.6.)

\proclaim{Theorem 8.1}
Suppose that $f$ is a proper holomorphic mapping of $\O$ onto the unit disc
with zeroes at $a_1,\dots,a_N$ with multiplicities $M(1),\dots,M(N)$,
respectively.  The weighted Szeg\H o kernel $\sigma(z,w)$ satisfies
$$\sigma(z,w)=\frac{1}{1-f(z)\overline{f(w)}}\left(
\sum_{i,j=1}^{N}
\sum_{n=0}^{M(i)}
\sum_{m=0}^{M(j)}
c_{ijnm}\sigma_{\bar n}(z,a_i)\,\overline{\sigma_{\bar m}(w,a_j)}\right).
\tag8.2$$
\endproclaim

The coefficients $c_{ijnm}$ can easily be determined.  Suppose
$1\le k\le N$ and $0\le q \le M(k)$.  Differentiate (8.2) $q$ times
with respect to $\bar w$ and set $w=a_k$ (and recall that $f$ has a
zero of multiplicity $M(k)$ at $a_k$) to obtain
$$\sigma_{\bar m}(z,a_k)=
\sum_{i,j=1}^{N}
\sum_{n=0}^{M(i)}
\sum_{m=0}^{M(j)}
c_{ijnm}\sigma_{\bar n}(z,a_i)\,\overline{\sigma_{q\bar m}(a_k,a_j)}$$
We saw an identity like this when we showed above that the functions $h_{jnp}$
are linearly independent for each fixed~$k$.  The same reasoning we used there
yields that such a relation can only be true if the system,
$$
\sum_{j=1}^{N}
\sum_{m=0}^{M(j)}
c_{ijnm}\overline{\sigma_{q\bar m}(a_k,a_j)}
=
\cases
1, &\text{if }i=k\text{ and } m=q,\\
0, &\text{if }i\ne k\text{ or } m\ne q,
\endcases
$$
has a unique solution.  This gives us a non-degenerate linear system 
to solve for the coefficients $c_{ijnm}$.

In the case that the proper holomorphic mapping $f$ in Theorem~8.1 has
simple zeroes, the formula for the weighted Szeg\H o kernel becomes
easier to write.

\proclaim{Theorem 8.2}
Suppose that $f$ is a proper holomorphic mapping of $\O$ onto the unit
disc with simple zeroes at $a_1,\dots,a_N$.
The weighted Szeg\H o kernel $\sigma(z,w)$ satisfies
$$\sigma(z,w)=\frac{1}{1-f(z)\overline{f(w)}}
\sum_{i,j=1}^{N}
c_{ij}\sigma(z,a_i)\,\overline{\sigma(w,a_j)}
\tag8.3$$
where the coefficients $c_{ij}$ are determined by the condition that
the matrix formed by the coefficients $[c_{ij}]$ is the inverse  
to the matrix $[\sigma(a_k,a_j)]$.
\endproclaim

We remark that for any proper holomorphic mapping $f$ from $\O$ onto
the unit disc, it is always possible to choose a M\"obius transformation
$\phi$ so that the proper map $\phi\circ f$ has simple zeroes.  Hence
the simpler formula in Theorem~8.2 is always at our disposal.

We can use identity (8.3) and (7.1) to derive similar results for
the weighted Garabedian kernel.  Assume that $w$ is in $b\O$ and
multiply (8.3) by $i\phi(w)\overline{T(w)}$.  Use (7.1)
to obtain
$$\lambda(w,z)=\frac{1}{1-f(z)\overline{f(w)}}
\sum_{i,j=1}^{N}
c_{ij}\sigma(z,a_i)\,{\lambda(w,a_j)}.$$
Finally, replace $\overline{f(w)}$ by $1/f(w)$ to obtain
$$\lambda(w,z)=\frac{f(w)}{f(w)-f(z)}
\sum_{i,j=1}^{N}
c_{ij}\sigma(z,a_i)\,{\lambda(w,a_j)}.$$
Since both sides of this identity are holomorphic in $w$, the identity
extends to hold for all $w$ in $\O$ minus the finite set of points
$\{a_1,\dots,a_N\}$ and $z$.  The singularities at $\{a_1,\dots,a_N\}$
are easily seen to be removable because $f$ vanishes at these points.

\proclaim{Theorem 8.3}
Suppose that $f$ is a proper holomorphic mapping of $\O$ onto the unit
disc with simple zeroes at $a_1,\dots,a_N$.
The weighted Garabedian kernel $\lambda(z,w)$ satisfies
$$\lambda(w,z)=\frac{f(w)}{f(w)-f(z)}
\sum_{i,j=1}^{N}
c_{ij}\sigma(z,a_i)\,{\lambda(w,a_j)}
\tag8.4$$
where the coefficients $c_{ij}$ are determined by the condition that
the matrix formed by the coefficients $[c_{ij}]$ is the inverse
to the matrix $[\sigma(a_k,a_j)]$.
\endproclaim

\subhead 9. The weighted Szeg\H o kernel and the double of a domain \endsubhead
Suppose that $\O$ is an $n$-connected domain in the plane such that no
boundary component is a point and assume further that the boundary of $\O$
consists of $n$ non-intersecting $C^\infty$ smooth
closed curves.  As before, we let $\Oh$ denote the double of $\O$ and
$R(z)$ denote the antiholomorphic involution on $\Oh$ which fixes the boundary
of $\O$ and let $\Ot=R(\O)$ denote the reflection of $\O$ in $\Oh$ across the
boundary.  Recall that if $f:\O\to D_1(0)$ is a proper holomorphic mapping
of $\O$ onto the unit disc, then $f$ extends to be a
meromorphic function on $\Oh$.

We shall now prove that for fixed points $A_1$ and $A_0$ in $\O$, functions
of $z$ of the form $\sigma(z,A_1)/\sigma(z,A_0)$ extend as meromorphic
functions to the double of $\O$.  Indeed, if we write the conjugate
of formula (7.1), first using $a=A_1$ and then $a=A_0$, and divide
the two resulting formulas, we see that
$\sigma(z,A_1)/\sigma(z,A_0)$ is equal to the complex conjugate of
$\lambda(z,A_1)/\lambda(z,A_0)$ when $z\in b\O$.  The arguments in
\cite{13} can easily be adapted to show that functions of $z$ of the form
$\sigma(z,a)$ cannot vanish to infinite order at a boundary point for
fixed points $a$ in $\O$.  Hence, the zeroes of
$\sigma(z,A_1)$ and $\sigma(z,A_0)$ on $b\O$ are isolated and of finite
order.  Now the usual reflection argument (which maps a one sided
neighborhood of a boundary point of $\O$ onto the lower half disc)
in the construction of the double of $\O$ shows that
$\sigma(z,A_1)/\sigma(z,A_0)$ is a meromorphic function on $\O$ which extends
up to $b\O$ with at most finitely many pole-like singularities on $b\O$,
and the conjugate of $\lambda(R(z),A_1)/\lambda(R(z),A_0)$
is a meromorphic function on $\Ot$ which extends up to $b\O$
from the ``outside'' of $\O$ and which agrees with
$\sigma(z,A_1)/\sigma(z,A_0)$ on $b\O$.  Hence, $\sigma(z,A_1)/\sigma(z,A_0)$
extends to $\Oh$ as a meromorphic function.
Similar reasoning using (7.1) and (7.2) shows that for fixed points
$A_1$ and $A_0$ in $\O$, functions of $z$ of the form
$\sigma_{\bar n}(z,A_1)/\sigma(z,A_0)$
extend as meromorphic functions to the double of $\O$.

If we fix a point $b$ in $\O$ and if we divide (8.2) by
$\sigma(z,b)\sigma(b,w)$, we see that
$$\frac{\sigma(z,w)}{\sigma(z,b)\sigma(b,w)}$$
is a rational combination of holomorphic functions of $z$ that extend
meromorphically to the double of $\O$ and antiholomorphic functions
of $w$ that extend antimeromorphically to the double of $\O$.  Since
the field of meromorphic functions on the double of $\O$ is
generated by just two such functions (a primitive pair), we obtain
the following result about the complexity of $\sigma(z,w)$.

\proclaim{Theorem 9.1}
Suppose that $\O$ is an $n$-connected domain in the plane such
that the boundary of $\O$ consists of $n$ non-intersecting
analytic $C^\infty$ smooth closed curves.  Let $G_1$ and $G_2$
be a primitive pair for the field of meromorphic functions on
the double of $\O$.  Given any point $b$ in $\O$,
the weighted Szeg\H o kernel is given by
$$\sigma(z,w)=\sigma(z,b)\sigma(b,w)
R(G_1(z),G_2(z),\overline{G_1(w)},\overline{G_2(w)}),$$
and thus is a rational combination of only three functions of
one complex variable on $\O$.
\endproclaim

We showed in \S8 that it is not possible for a function of the form
$$\sigma(z,A_1)-c\sigma(z,A_0)$$
to vanish identically on $\O$.  (If it did, then every polynomial $p(z)$ would
satisfy $p(A_1)-\bar c p(A_0)=0$ by virtue of the reproducing property of the
weighted Szeg\H o kernel, and this is absurd.)  Hence, if $A_1\ne A_0$, then
$\sigma(z,A_1)/\sigma(z,A_0)$ extends to $\Oh$ as a non-constant meromorphic
function of some finite order $m$ on $\Oh$.  We record this result
here for future use.

\proclaim{Theorem 9.2}
Suppose that $\O$ is an $n$-connected domain in the plane such
that the boundary of $\O$ consists of $n$ non-intersecting
$C^\infty$ smooth closed curves.  If $A_1$ and $A_0$ are distinct points in
$\O$, then $\sigma(z,A_1)/\sigma(z,A_0)$ extends to the double of $\O$ as a
non-constant meromorphic function.
\endproclaim

\subhead 10. A special weight function for the Szeg\H o kernel\endsubhead
Assume, as we did in \S9, that $\O$ is an $n$-connected domain in the
plane such that the boundary of $\O$ consists of $n$ non-intersecting
$C^\infty$ smooth closed curves.  Let $p(a,z)$ denote
the classical Poisson kernel associated to $\O$ which reproduces
harmonic functions in the sense that
$$u(a)=\int_{z\in b\O}p(a,z)\,u(z)\ ds$$
when $u$ is harmonic in $\O$ and continuous up to the boundary.
Choose a point $A_0$ in $\O$ and define a weight function via
$$\phi(z)=p(A_0,z).$$
The weighted Szeg\H o kernel associated to this weight has the virtue
that
$$\sigma(z,A_0)\equiv 1.$$
Formula (7.1) shows that $\lambda(z,A_0)$ is non-vanishing on the
boundary of $\O$ and the argument principle and formula (7.1) show
that $\lambda(z,A_0)$ has exactly $n-1$ zeroes in $\O$ (counted with
multiplity) in the $z$ variable.  Furthermore, since
$$p(w,z)=\frac{1}{2\pi}\frac{\dee}{\dee n_z}G(w,z)=
\frac{i}{\pi}\frac{\dee}{\dee\bar z}G(w,z)\overline{T(z)},$$
it follows that
$$p(A_0,z)=
\frac{i}{\pi}\frac{\dee}{\dee\bar z}G(A_0,z)\overline{T(z)},$$
and (7.1) reveals that
$$\lambda(z,A_0)=
\frac{1}{\pi}\frac{\dee}{\dee\bar z}G(A_0,z)\tag 10.1$$
for $z\in b\O$.  Since these two functions extend meromorphically
inside $\O$ and have exactly the one simple pole at $A_0$ with exactly
the same residue, they are equal for all $z$ in $\O$, too.

Because $\sigma$ can be expressed as
the weighted projection of the weighted Cauchy kernel (see \S7),
$\sigma(z,a)=P_\phi C_a$, and because $P_\phi$ is a continuous linear
operator from $C^\infty(b\O)$ into $A^\infty(\O)$, we may conclude
that $\sigma(z,a)$ has no zeroes in the $z$ variable in $\Obar$ for
all $a$ sufficiently close to $A_0$.  Furthermore, $\lambda(z,a)$ is
non-vanishing on $b\O$ and has exactly $n-1$ zeroes in $\O$ when $a$
is close to $A_0$.  Choose a point $A_1\ne A_0$ which is close enough
so that these two conditions hold.  Let $G(z)$ denote the extension of
$$\sigma(z,A_1)=\frac{\sigma(z,A_1)}{\sigma(z,A_0)}$$
to the double of $\O$ given by Theorem~9.2.  Notice that
$$G(z)=\overline{\lambda(R(z),A_1)}/\overline{\lambda(R(z),A_0)}$$
on $\Ot$.  Hence, $G$ has no zeroes and no poles in $\Obar$ and
some positive number $m\le n-1$ of zeroes in $\Ot$ and the same number $m$
of poles in $\Ot$.  The order of $G$ on $\Oh$ is $m$.

Choose a point $w_0\ne0$ in $\C$  close enough to the origin that
$G^{-1}(w_0)$ consists of $m$ distinct points in $\Oh$ which all
fall in $\Ot$.  We may also choose $w_0$ so that none of these points
is one of the $n-1$ zeroes of $\lambda(R(z),A_0)$ in $\Ot$.  Let
$z_1,\dots,z_m$ denote these $m$ points.

We now wish to show that it is possible to choose a point $A_2$ in
$\O$ so that the meromorphic extensions of $\sigma(z,A_1)$ and
$$\sigma(z,A_2)=\frac{\sigma(z,A_2)}{\sigma(z,A_0)}$$
to $\Oh$ form a primitive pair for $\Oh$ (meaning that they generate
the field of meromorphic functions on the double of $\O$).
For a fixed point $A_2$ in $\O$, let $H(z)$ be
defined to be the meromorphic extension of $\sigma(z,A_2)$
to $\Oh$.  Think of $H$ as depending on $A_2$ even though we have suppressed
this fact in the notation.  To show that $G(z)$ and $H(z)$
form a primitive pair, we need only choose $A_2$ so that $H(z)$ separates
the $m$ points in  $G^{-1}(w_0)$ (see \cite{3, page~321-324}).

Suppose that $A_2$ is not equal to $A_1$ and suppose that $A_2$ is close
enough to $A_0$ that $\lambda(R(z_k),A_2)\ne0$ for each $k$.

If $z_i\in\Ot$ and $z_j\in\Ot$ are points in $G^{-1}(w_0)$ which are not
separated by $H$, then 
$$\frac{\lambda(R(z_i),A_2)}{\lambda(R(z_i),A_0)}=
\frac{\lambda(R(z_j),A_2)}{\lambda(R(z_j),A_0)},$$
and so
$$\frac{\lambda(R(z_i),A_2)}{\lambda(R(z_j),A_2)}=c$$
where
$$c=\frac{\lambda(R(z_i),A_0)}{\lambda(R(z_j),A_0)}$$
is a non-zero constant.
But the set of points $w$ in $\O$ where
$$\frac{\lambda(R(z_i),w)}{\lambda(R(z_j),w)}=c$$
is a finite subset of $\O$ because this function of $w$ extends to the
double of $\O$ as a non-constant meromorphic function since
$$\frac{\lambda(R(z_i),w)}{\lambda(R(z_j),w)}$$
is equal to the conjugate of
$$\frac{\sigma(R(z_i),R(w))}{\sigma(R(z_j),R(w))}$$
on $b\O$.  Hence $w=A_2$ can be chosen to avoid this
possibility for each pair of indices $i\ne j$.

Assume that $A_1$ and $A_2$ in $\O$ are two points in $\O$
such that the extensions of $\sigma(z,A_1)$ and $\sigma(z,A_2)$ to
the double of $\O$ form a primitive pair for the field of meromorphic
functions on the double.

Let $f$ denote an Ahlfors mapping of $\O$ onto the unit disc.
We may suppose that the base
point of the map has been chosen so that the zeroes of $f$ are all
simple zeroes.  We know that $f(z)$ extends meromorphically
to $\Oh$.  It follows that all the functions that appear on the right
hand side of formula (8.3) extend to the double of $\O$.  Hence,
$\sigma(z,w)$ is a rational combination of $\sigma(z,A_1)$, $\sigma(z,A_2)$
and the conjugates of $\sigma(w,A_1)$ and $\sigma(w,A_2)$.

We collect these results in the following theorem.

\proclaim{Theorem 10.1}
Suppose that $\O$ is an $n$-connected domain in the
plane such that the boundary of $\O$ consists of $n$ non-intersecting
$C^\infty$ smooth closed curves.
The weighted Szeg\H o kernel with respect to the Poisson weight
$$\phi(z)=p(A_0,z)$$
for a point $A_0$ in $\O$ is such that $\sigma(z,a)$ extends to the
double of $\O$ as a non-constant meromorphic function for each $a$ in
$\O$.  Furthermore, there exist two points $A_1$ and $A_2$ in $\O$
such that the extensions of $\sigma(z,A_1)$ and $\sigma(z,A_2)$ to
the double of $\O$ form a primitive pair for the field of meromorphic
functions on the double.  The kernel $\sigma(z,w)$ is a rational
function of
$\sigma(z,A_1)$, $\sigma(z,A_2)$,
$\overline{\sigma(w,A_1)}$ and $\overline{\sigma(w,A_2)}$.
\endproclaim

\subhead 11. Finite Riemann surfaces\endsubhead
Suppose that $\O$ is a finite Riemann surface with boundary, i.e.
suppose that $\O$ is a Riemann surface with boundary with compact
closure of finite genus and finitely many boundary curves.  We assume
that there is at most one boundary curve and that none of the boundary
curves are pointlike.  To make the exposition easier, we shall assume
that the boundary curves of $\O$ are $C^\infty$ smooth real analytic
curves.  It will be clear that this assumption can be greatly reduced
in what follows, but we do not concern ourselves with this
here.  It is a standard construction to produce the Szeg\H o projection
and kernel with respect to a weight function on the boundary of
$\O$ (see \cite{16}).  Given a measure $\omega$ on the boundary which is
given by a positive $C^\infty$ weight function with respect to the
standard metric on the boundary, let $S(z,w)$ denote the Szeg\H o kernel
function defined on $\O\times b\O$ which reproduces holomorphic functions
on $\O$ with respect to $\omega$ in the sense that
$$h(z)=\int_{w\in b\O}S(z,w)h(w)\,d\omega$$
for points $z$ in $\O$ and holomorphic functions $h$ on
$\O$ that are in the $L^2$ Hardy space associated to $\O$ relative to
the measure $\omega$.  Let $P$ denote the Szeg\H o projection, which
is the orthogonal projection of the $L^2$ space on $b\O$ with respect
to the weight function $\omega$ onto the closed subspace of functions
in this space which are the boundary values of holomorphic functions
on $\O$.

Ahlfors proved that Ahlfors maps exist in the more general setting
that we are dealing with now (see \cite{2}).  Fix a point $a$ in $\O$
and let $f_a$ denote an Ahlfors map associated to $(\O,a)$.  This map
is a holomorphic map of $\O$ onto the unit disc such that $f_a(a)=0$
and which maximizes the derivative $f'_a(a)$ in some coordinate chart.
Ahlfors proved that this map is a proper holomorphic mapping of $\O$
onto the unit disc.  Since the boundary of $\O$ is $C^\infty$ smooth,
it follows that $f_a$ extends $C^\infty$ smoothly up to the boundary of $\O$.
Of course, $f_a$ maps the boundary of $\O$ into the boundary of the
unit disc, i.e. $|f(z)|=1$ when $z\in b\O$.

Exactly the same arguments as those given in \S\S7-8 can now be applied in
this more general context to yield that the Szeg\H o kernel $S(z,w)$ is a
rational combination of finitely many holomorphic functions of one
complex variable on $\O$.  Indeed, if the zeroes of $f_a$ are simple,
then the formula in Theorem~8.2 holds with $S(z,w)$ in place of
$\sigma(z,w)$ and $f_a$ in place of $f$.  If the zeroes are not simple,
then the formula in Theorem~8.1 holds where it is understood that
the derivatives $S_{\bar n}(z,a_j)$ in the second variable are taken
with respect to some arbitrary, but fixed, coordinate chart near $a_j$.

To show that $S(z,w)$ is actually generated by {\it only three\/} holomorphic
functions of one variable, we must do a little extra work.  Let
$G(z,w)$ denote the classical Green's function for $\O$.  For a fixed
point $b$ in $\O$, let $\dee_z G(z,b)$ denote the meromorphic one-form
$(\dee/\dee z)G(z,b)\,dz$.  This form reproduces holomorphic functions
on $\O$ in the sense that
$$h(b)=\int_{z\in b\O}h(z) \dee_z G(z,b).$$
There is a $C^\infty$ function ${\Cal G}_b$ on the boundary of $\O$,
such that
$$\int_{z\in b\O}\psi(z) \dee_z G(z,b)=
\int_{z\in b\O}\psi(z)\overline{{\Cal G}_b(z)}\,d\omega$$
for all continuous functions $\psi$ on $b\O$.
It follows that the function $S_b(z):=S(z,b)$ is the Szeg\H o
projection of ${\Cal G}_b$.  Since $(S_b(z)-{\Cal G}_b)\omega$ is
orthogonal to holomorphic functions, we may use a
theorem of Read (see \cite{16, page~75} and \cite{17}) to assert that
there is a meromorphic one-form $A_b$ on $\O$ with no singularities on
$\Obar$ such that
$$\int_{z\in b\O}\psi(z)\left(S_b(z)-\overline{{\Cal G}_b(z)}\right)\,d\omega=
\int_{z\in b\O}\psi(z)*A_b$$
for all continuous functions $\psi$ on $b\O$.
Since $S_b$ and ${\Cal G}_b$ are $C^\infty$ smooth on $b\O$, it
follows that $A_b$ is $C^\infty$ smooth up to $b\O$ and hence,
there is a $C^\infty$ smooth function $\alpha_b(z)$ on $b\O$
such that
$$\int_{z\in b\O}\psi(z)*A_b=
\int_{z\in b\O}\psi(z)\overline{\alpha_b(z)}\,d\omega$$
for all continuous functions $\psi$ on $b\O$.  Define $\lambda_b(z)$ to be
${\Cal G}_b(z)+\alpha_b(z)$.  For two points $a$ and $b$ in $\O$,
the quotient $S_b(z)/S_a(z)$ is a meromorphic function on $\Obar$.
It is easy to verify that the quotient $\lambda_b(z)/\lambda_a(z)$
is the boundary value of an antimeromorphic function on $\Obar$ given
as the quotient of two meromorphic one-forms.  Hence
$S_b(z)/S_a(z)$ extends to the double of $\O$ as a meromorphic function.
This and the fact that we may compose by a M\"obius transformation so
as to be able to assume that our proper holomorphic map to the disc
has simple zeroes is all that we need to be able to use the formula in
Theorem~8.2 to deduce that
$$S(z,w)=S(z,a)S(a,w)R(G_1(z),G_2(z),\overline{G_1(w)},\overline{G_2(w)})$$
where $R$ is a rational function and $G_1$ and $G_2$ form a primitive
pair for the double of $\O$.

The arguments in \S10 carry over to our finite Riemann surface and
Theorem~10.1 holds in this more general context.

We conclude by showing how these results can be applied to the Bergman
kernel on $\O$.  The Bergman kernel on $\O$ is a differential $(1,1)$
form given by
$$K(w,z)dw\wedge d\bar z=
\frac{\dee^2}{\dee w\dee\bar z}G(w,z)dw\wedge d\bar z.$$
Let $\alpha=df/f$ where $f$ is a proper holomorphic mapping of $\O$ onto
the unit disc (such as an Ahlfors map).  The proof of Theorem~1.2 given in
\cite{9} yields that
$$\frac{K(w,z)dw\wedge d\bar z}{\alpha(w)\wedge\overline{\alpha(z)}}$$
can be viewed as a {\it function\/} on $\O\times\O$ which extends to
$\Oh\times\Oh$ to be meromorphic in $w$ and antimeromorphic in $z$.
It now follows from the generalization of Theorem~10.1 mentioned above that
the Bergman kernel associated to $\O$ is given as
$$K(w,z)dw\wedge d\bar z=
R(G_1(w),G_2(w),\overline{G_1(z)},\overline{G_2(z)})
df(w)\wedge\overline{df(z)}$$
where $R$ is a complex rational function and $G_1$ and $G_2$ form
a primitive pair for $\Oh$.  It is interesting to note
that the two functions $G_1(z)$ and $G_2(z)$ can be taken to be
$S(z,A_1)$ and $S(z,A_2)$ for suitably chosen points $A_1$ and $A_2$ in
$\O$ and that $f$ is also a rational combination of these two functions.

I leave it for the future to relate other objects of potential theory
associated to a finite Riemann surface to this weighted Szeg\H o kernel
so as to obtain results about complexity in complex analysis and
potential theory.

\enddocument